\documentclass[a4paper]{amsart}

\usepackage{a4}
\usepackage[utf8]{inputenc}
\usepackage[T1]{fontenc}
\usepackage{lmodern}
\usepackage{amsmath,amssymb,amsfonts}
\usepackage{amsthm}
\usepackage{mathtools}
\usepackage{mathrsfs}
\usepackage{graphicx}
\usepackage{url}
\usepackage{xcolor}
\usepackage{textcomp}
\usepackage{enumerate}
\usepackage{multirow}
\usepackage{subcaption}
\usepackage{makecell}
\usepackage{booktabs}
\usepackage{algorithm}
\usepackage{algorithmicx}
\usepackage{algpseudocode}
\usepackage{listings}
\usepackage{scalefnt}
\usepackage{placeins}

\newcommand{\x}{\mathrm{x}}
\newcommand{\lossv}{\mathrm{Loss}_{\tiny\mbox{$V$}}}
\newcommand{\lossg}{\mathrm{Loss}_{\tiny\mbox{$\Phi$}}}
\newcommand{\lossgradv}{\mathrm{Loss}_{\tiny\mbox{$\nabla_{x}V$}}}
\newcommand{\lossHJB}{\mathrm{Loss}_{\scalebox{0.5}{\text{HJB}}}}
\newcommand{\lossODE}{\mathrm{Loss}_{\scalebox{0.5}{\text{ODE}}}}

\def\dd{{\rm d}}

\def\B{\mathcal{B}}

\def\H{\mathcal{H}}

\def\P{\mathcal{P}}

\def\T{\mathcal{T}}

\def\W{\mathcal{W}}

\newcommand{\NN}{\mathbb{N}}

\newcommand{\RR}{\mathbb{R}}

\newcommand\be{\begin{equation}}
	\newcommand\ee{\end{equation}}
\newcommand\ba{\begin{array}}
	\newcommand\ea{\end{array}}
\newcommand{\bean}{\begin{eqnarray*}}
	\newcommand{\eean}{\end{eqnarray*}}

\def\ds{\displaystyle}

\title[High-dimensional optimal control and mean field games]{Initialization-driven neural generation and training for high-dimensional optimal control and first order mean field games}

\author{Mouhcine Assouli}
\address{Moroccan Center for Game Theory, UM6P, Rocade, Rabat-Sal\'e, 11103, Morocco}
\address{Institut de recherche XLIM-DMI, UMR 7252 CNRS, Facult\'e des Sciences et Techniques, Universit\'e de Limoges, Limoges, 87060, France}
\email{mouhcine.assouli@um6p.ma}

\author{Justina Gianatti}
\address{CIFASIS (CONICET-UNR), Bv. 27 de Febrero 210 bis, Rosario, S2000EZP, Argentina}
\email{gianatti@cifasis-conicet.gov.ar}

\author{Badr Missaoui}
\address{Moroccan Center for Game Theory, UM6P, Rocade, Rabat-Sal\'e, 11103, Morocco}
\email{badr.missaoui@um6p.ma}

\author{Francisco J. Silva}
\address{Institut de recherche XLIM-DMI, UMR 7252 CNRS, Facult\'e des Sciences et Techniques, Universit\'e de Limoges, Limoges, 87060, France}
\email{francisco.silva@unilim.fr}


\begin{document}

\begin{abstract}
We introduce a neural-network-based method for approximating the value function of high-dimensional deterministic optimal control problems. The proposed approach exploits the relationship between Pontryagin's Maximum Principle (PMP) and the value function, which is characterized as the unique viscosity solution of the Hamilton–Jacobi–Bellman (HJB) equation. A neural network is first trained to obtain a coarse approximation of the value function by minimizing the residual of the HJB equation. The gradient of this approximation is then used to initialize the numerical solution of the two-point boundary value problem arising from PMP, enabling the generation of reliable optimal trajectories, adjoint states, and costs. This dataset is then used to further train the neural network through a loss function enforcing the HJB equation.

We next address the computation of equilibria in first-order Mean Field Game (MFG) problems by integrating the proposed methodology with the fictitious play algorithm. Such equilibria are characterized by a coupled system consisting of an HJB equation and a continuity equation. To approximate the solution of the continuity equation, we introduce a second neural network that learns the flow map transporting the initial agent distribution along optimal trajectories. This network is trained using data obtained by solving the ordinary differential equations (ODEs) associated with the controlled dynamics. The iterative procedure is initialized via joint training of the value function and the flow map, based on a composite loss involving both HJB and ODE residuals. Several numerical experiments demonstrate the accuracy and robustness of the proposed method.
\end{abstract}

\maketitle

\section{Introduction}
\label{Sect_Intro} 

Machine learning methods for optimal control have attracted considerable
attention over the past decade. A central objective is the construction of
robust optimal feedback controllers, commonly approached by approximating the
solution to the associated Hamilton--Jacobi--Bellman (HJB)
equation~\cite{bardi1997optimal,falcone2013semi}. For finite-horizon
deterministic problems, the HJB equation is a first-order nonlinear PDE
describing the optimal cost as a function of the initial time and state.
Classical methods, including finite difference, semi-Lagrangian, and finite
element schemes, rely on spatial grids and therefore suffer from the
\emph{curse of dimensionality}~\cite{MR134403}. We refer
to~\cite{falcone2013semi,MR3653852} and the references therein for an overview
of grid-based methods for HJB equations. Several alternatives have been
proposed for high-dimensional problems, including tropical
methods~\cite{MR2346381,MR2385864,MR2599910,MR4546175}, semi-Lagrangian schemes
on sparse grids~\cite{MR3045704,MR3592131}, polynomial
approximations~\cite{kang2017mitigating,MR3769705,MR4109008,MR4253741,MR4664738},
optimization methods based on the Hopf and Lax--Oleinik
formulae~\cite{MR3413587,MR3543239,MR4221591,MR4754321,MR4717768,MR4713729},
tree-structured semi-Lagrangian schemes~\cite{MR3984311,MR4506570}, tensor
decomposition
techniques~\cite{Horowitz14,Gorodetsky_2018,MR4254983,MR4442445,MR4729053,MR4641650},
and neural network (NN)
approximations~\cite{MR4123395,MR4081911,MR4338293,MR4233238,MR4166062,
nakamura2021adaptive,MR4275044,onken2022neural,MR4646437,MR4907583,
MR4793480,MR4834792}.

Efficient solvers for high-dimensional HJB equations also provide a natural
basis for tackling high-dimensional Mean Field Games (MFGs). Introduced
independently by J.-M.~Lasry and P.-L.~Lions
in~\cite{LasryLions06i,LasryLions06ii,MR2295621} and by Caines, Huang, and
Malham\'e in~\cite{HMC06}, MFGs describe the limiting behavior of Nash
equilibria in symmetric stochastic games as the number of players tends to
infinity. We refer to~\cite{MR3195844,MR3559742,MR3752669,MR3753660,MR4214773}
for an overview of the theory and its applications to crowd motion, economics,
and finance. In their simplest form, MFG equilibria are characterized by a
coupled system consisting of an HJB equation for the value function of a
representative player and a transport equation for the distribution of the
agents. Numerical approaches include finite difference
schemes~\cite{AchdouCapuzzo10,MR3097034,MR3452251}, semi-Lagrangian
schemes~\cite{MR3148086,CarliniSilva18,MR3392626,Silva2024LG}, finite element
methods~\cite{MR4688199,MR4858134,MR4854642}, approximations by finite-state
discrete-time MFGs~\cite{MR4030259,MR4835150,MR4714572}, and, for
high-dimensional problems, machine learning
techniques~\cite{MR4264647,MR4236167,MR4522347,MR4587578,assouli2023deep,
assouli2023policy}. Numerical methods for MFGs, including variational
approaches, are reviewed in~\cite{MR4214777,MR4368188}.

Our first contribution is the introduction of the
Initialization--Generation--Training (IGT) method. As
in~\cite{kang2015causality,kang2017mitigating,nakamura2021adaptive,MR4253741,
MR4442445}, it constructs an approximation of the value function from data
generated by solving a family of open-loop optimal control problems
parameterized by the initial time and state. Each problem is treated using
Pontryagin's Maximum Principle (PMP)~\cite{MR166037}, which yields a two-point
boundary value problem (TPBVP) providing necessary conditions for local
optimality. The main distinction from previous approaches is the systematic
use of sensitivity relations from optimal control theory. When the value
function is sufficiently smooth, these relations imply that a solution to the
TPBVP is globally optimal if its adjoint state coincides, along the state
trajectory, with the spatial gradient of the value function; see
Section~\ref{sec:preliminaries}. This suggests using an approximation of this
gradient to initialize the TPBVP solver and thereby favor convergence toward
globally optimal solutions.

A related strategy was proposed in~\cite{cristiani2010initialization}, where a finite difference approximation of the value function is used to initialize a shooting method. Its reliance on spatial grids, however, restricts the approach to low-dimensional problems. The IGT framework instead uses a mesh-free NN approximation. More precisely, we first compute a rough approximation of the value function using the Deep Galerkin Method (DGM)~\cite{sirignano2018dgm}, by minimizing the residual of the HJB equation. The spatial gradient of this approximation is then used to provide a suitable initialization of Newton's iterates to approximate global solutions to the parameterized family of optimal control problems. This procedure generates reliable data including optimal costs, optimal states, and adjoint states, which are used to retrain the NN. The improved approximation yields feedback controls and can be used to generate additional data in subsequent IGT rounds.

The Adaptive Sampling and Model Refinement (ASMR) algorithm introduced
in~\cite{nakamura2021adaptive}, and inspired by~\cite{MR4275044}, also seeks to
reduce sensitivity to the TPBVP initialization by using a time-marching
procedure through intermediate times. However, this procedure is not
specifically designed to identify globally optimal solutions and fails to
converge in some of the experiments reported in
Section~\ref{sec:numerical_results}. We provide extensive numerical evidence
of the effectiveness of IGT by comparing the learned value functions and
feedback controls with explicit solutions. The test cases include smooth and
non-smooth problems in several state-space dimensions. Overall, IGT accurately
recovers both the value function and the optimal feedback, whereas alternative
approaches may be less accurate or fail as the dimension increases or the
value function loses regularity.

The benefit of the proposed initialization is particularly clear in the
obstacle-avoidance experiment of Section~\ref{subsec:obstacle_oc}. During the
first IGT round, the TPBVP solver converges for $14$, $8$, and $3$ of the $32$
initial states considered in dimensions $d=2$, $10$, and $50$, respectively.
After a few rounds, it converges for all $32$ initial states in every
dimension. By contrast, arbitrary initialization fails for all tested initial
states, even when combined with time-marching. These results highlight the
decisive role of the initialization in the convergence of the TPBVP solver.

Our second contribution extends IGT to the approximation of first-order MFG
equilibria by coupling it with an NN approximation of solutions to continuity
equations. For deterministic games, the system introduced
in~\cite{MR2295621} consists of a first-order HJB equation coupled with a
continuity equation. Under suitable structural assumptions, the fictitious
play method developed for MFGs in~\cite{MR3608094,These_Saeed_18} approximates
an equilibrium by solving these two equations separately and iteratively. In
our framework, the HJB equation is treated with IGT, while the continuity
equation is handled through a learned flow map.

Indeed, the solution to the continuity equation is the push-forward of the
initial distribution by the flow generated by the feedback control associated
with the HJB equation. We therefore approximate this flow by an NN depending
on the current approximation of the value function. At each iteration,
training data are generated by solving the underlying controlled ODEs for
initial states sampled from the initial distribution. The learned flow
provides a compact generator of the population distribution, avoiding the
need to store the growing collection of trajectories produced by the
fictitious play iterations. The procedure is initialized by jointly training
the value function and the flow map through a composite loss involving the HJB
and ODE residuals. In the numerical experiments, convergence is monitored
using both the aggregated Sinkhorn divergence~\cite{genevay2018learning}
between the current distribution and its best response, and the
\emph{exploitability}~\cite{Lauriere_fictitious_play}.

The paper is organized as follows. Section~\ref{sec:preliminaries} introduces
the optimal control and MFG problems considered throughout the paper and
recalls the required preliminary results. Section~\ref{sec:IGT} presents the
IGT method for finite-horizon deterministic optimal control problems.
Section~\ref{sec:MFG_approximation} combines IGT with the NN approximation of
continuity equations and fictitious play. Finally,
Section~\ref{sec:numerical_results} reports numerical experiments for
high-dimensional optimal control problems, including linear-quadratic and
non-smooth examples, quadcopter guidance, and obstacle avoidance. The
initialization step is particularly important in the last two problems. The
section concludes with applications to a linear-quadratic MFG with an explicit
solution and to an MFG involving obstacle avoidance and congestion aversion.
	
\section{Preliminaries on deterministic optimal control problems and MFG systems}
\label{sec:preliminaries}
Given $T>0$, a nonempty closed subset $A$ of $\RR^{m}$, $\ell\colon[0,T]\times\RR^{d}\times A\to\RR$, $g\colon\RR^{d}\to\RR$, and $b\colon[0,T]\times\RR^{d}\times A\to\RR^{d}$, we consider the following family of optimal control problems, parameterized by the initial time $t\in [0,T]$ and the initial state $x\in\RR^{d}$: 
	\begin{equation}
		\tag{\mbox{${\bf P}^{t,x}$}}
		\label{def:oc_problem}
		\begin{aligned}
			\inf\ & \int_{t}^{T} \ell(s, \mathrm{x}(s), \alpha(s))\, \mathrm{d}s + g\left( \mathrm{x}(T) \right) \\
			\text{s.t.} \quad & \dot{\mathrm{x}}(s) = b(s, \mathrm{x}(s), \alpha(s)) \quad \text{for a.e. } s \in [t,T], \\
			& \mathrm{x}(t) = x, \\
			& \alpha(s) \in A \quad \text{for a.e. } s \in [t,T].
		\end{aligned}
	\end{equation}
	Under standard assumptions on the data $(\ell,g,b)$, problem~\eqref{def:oc_problem} is well-defined and its optimal value is finite (see, e.g.,~\cite[Section I.3]{fleming2006controlled}).  A key result in optimal control theory is Pontryagin's maximum principle (PMP)~\cite{MR166037}, which states that, under some differentiability assumptions over $(\ell,g,b)$ (see e.g.,~\cite[Section I.6]{fleming2006controlled}), for every solution $(\x^*,\alpha^*)$ to~\eqref{def:oc_problem} there exists $p^*\colon [t,T]\to\RR^{d}$, called {\it adjoint state}, such that 
	\begin{equation}
		\left\{\begin{aligned}
			\dot{\x}^*(s) &= b\left(s, \x^*(s), \alpha^*(s)\right)\quad\text{for a.e. }s\in [t,T],\\
			\x^*(t)&=x, \\
			\dot{p}^*(s) &= \nabla_{x}H\left(s, \x^*(s), p^*(s), \alpha^*(s)\right)\quad\text{for a.e. }s\in [t,T],\\
			p^*(T)&=\nabla g(\x^*(T)), \\
			\alpha^*(s) &\in \text{argmax}_{\alpha\in A} H(s, \x^*(s), p^*(s),\alpha)\quad\text{for a.e. }s\in [t,T],
		\end{aligned}
		\right.
		\label{eq:pmp_avec_H}  
	\end{equation}
	where the {\it pseudo-Hamiltonian} $H\colon[0,T]\times\RR^{d}\times\RR^{d}\times A\to\RR$ is given by 
	$$
	H(t,x,p,\alpha)=-\ell(t,x,\alpha)-p\cdot b(t,x,\alpha)\quad\text{for all }(t,x,p,\alpha)\in [0,T]\times\RR^{d}\times\RR^{d}\times A
	$$
	and $\nabla_{x}H$ denotes the gradient of $H$ with respect to the variable $x$. 
	
	In what follows, we suppose that, for every $(t,x,p)\in[0,T]\times\RR^{d}\times\RR^{d}$, the function $A\ni\alpha\mapsto H(t,x,p,\alpha)\in\RR$ admits a unique maximizer $\Psi(t,x,p)$, i.e., 
	\begin{equation}
		\{\Psi(t,x,p)\}=\text{argmax}_{\alpha\in A} H(t,x,p,\alpha).
		\label{eq:unique_minimizer_hamiltonian}
	\end{equation}
	
	Under this assumption,~\eqref{eq:pmp_avec_H} reduces to the following TPBVP:
	\begin{equation}
		\left\{\begin{aligned}
			\dot{\x}^*(s) &= b\left(s, \x^*(s), \Psi(s,\x^*(s),p^{*}(s))\right)\quad\text{for a.e. }s\in [t,T],\\
			\x^*(t)&=x, \\
			\dot{p}^*(s) &= \nabla_{x}H\left(s, \x^*(s), p^*(s), \Psi(s,\x^*(s),p^{*}(s))\right)\quad\text{for a.e. }s\in [t,T],\\
			p^*(T)&=\nabla g(\x^*(T)).
		\end{aligned}
		\right.
		\label{eq:TPBVP}  
	\end{equation}
	
	System~\eqref{eq:TPBVP} is a necessary condition for optimality and hence, finding solutions $(\x^*,p^*)$ to this system may provide optimal open-loop controls through the formula $\alpha^*(t)=\Psi(t,\x^*(t),p^*(t))$ for all $t\in [0,T]$. This procedure is at the heart of Pontryagin's approach to solve~\eqref{def:oc_problem}.
	
	Denote by $V(t,x)$ the optimal value of~\eqref{def:oc_problem}. Another fundamental result in optimal control theory (see e.g.~\cite[Chapter IV]{Fleming_Rishel_75}), having its roots in the pioneering work by R. Bellman~\cite{Bellman_57} on the dynamic programming principle, states that if $V$ is smooth enough, then it is characterized as the unique solution to the following Hamilton-Jacobi-Bellman (HJB) equation: 
	\begin{equation}
		\begin{aligned}   
			-\partial_{t} V(t,x) +\H  (t,x,\nabla_{x}V(t,x))&= 0\quad\text{for }(t,x)\in ]0,T[\times \RR^d, \\
			V(T,x) &=g(x) \quad\text{for }x\in  \RR^d,
		\end{aligned}
		\label{eq:hjb}
	\end{equation}
	where the {\it Hamiltonian} $\H\colon[0,T]\times\RR^{d}\times\RR^{d}\to\RR$ is defined by 
	$$
	\H(t,x,p)=\max_{\alpha\in A}H(t,x,p,\alpha)\quad\text{for all }(t,x,p)\in ]0,T[\times\RR^{d}\times\RR^{d}. 
	$$
	
	Let us mention that if $V$ is nonsmooth, but continuous, important results due to M. G. Crandall and P.-L. Lions~\cite{Crandall_Lions_81,Crandall_Lions_83} show that, under rather general assumptions on the data, $V$ is still the unique solution to~\eqref{eq:hjb} but in the so-called {\it viscosity} sense. 
	
	Interestingly, Pontryagin's maximum principle and the HJB approaches are connected through what is known as {\it sensitivity} analysis and {\it verification} result. Assuming that $V$ is smooth enough, the former states that if $(\x^*,\alpha^*)$ solves~\eqref{def:oc_problem}, and $p^*$ denotes the associated adjoint state, then 
	\begin{equation}
		\nabla_{x}V(s,\x^*(s))=p^*(s)\quad\text{for all }s\in[t,T],
		\label{eq:sensitivity}    
	\end{equation}
	while the latter asserts that $(\x^*,\alpha^*)$ is optimal for~\eqref{def:oc_problem} if and only if, for a.e. $s\in[t,T]$,
	\begin{equation}
		\alpha^*(s)=\Psi\left(s,\x^*(s),\nabla_xV(s,\x^*(s))\right).
		\label{eq:alpha_Psi_optimal}
	\end{equation}
	Altogether, if $V$ is smooth enough, then $(\x^*,\alpha^*)$ is a  (global) optimal solution to~\eqref{def:oc_problem} if and only if~\eqref{eq:pmp_avec_H} and~\eqref{eq:sensitivity} hold. This observation is the main motivation of the numerical approximation that we consider in Section~\ref{sec:IGT} below when solving~\eqref{def:oc_problem} over a batch of initial times and initial states. 
	
	Having at our disposal a smooth value function $V$, it follows from~\eqref{eq:alpha_Psi_optimal} that one can construct optimal feedback laws $[0,T]\times\RR^{d}\ni(t,x)\mapsto\alpha^{*}(t,x)\in A$ through the expression
	\begin{equation}
		\alpha^{*}(t,x)=\Psi\left(t,x,\nabla_{x}V(t,x)\right)\quad\text{for all }(t,x)\in[0,T]\times\RR^{d}
		\label{eq:optimal_feedbacks}
	\end{equation}
	and that the optimal dynamics for \eqref{def:oc_problem} is given by 
	\begin{equation}
		\left\{\begin{aligned}
			\dot{\x}^*(s) &= b\left(s, \x^*(s), \Psi\left(s,\x^*(s),\nabla_{x}V(s,\x^*(s))\right)\right)\\
			&=-\partial_{p}\H\left(s,\x^*(s),\nabla_{x}V(s,\x^*(s))\right)\quad\text{for a.e. }s\in [t,T],\\
			\x^*(t)&=x.
		\end{aligned}
		\right.
		\label{eq:optimal_dynamics}
	\end{equation}
	We refer the reader to~\cite{clarke1987relationship,zhou1990maximum,Zhou_93_verification} (see also~\cite[Chapter 5]{Yong_Zhou_book} and the references therein) for extensions of the sensitivity analysis and verification type results when the value function $V$ is nonsmooth. 
	
	To approximate the value function $V$, we consider a neural network $V_{\theta}$ parameterized by $\theta$. The parameter $\theta$ is computed by minimizing a suitable loss function, as described in Section~\ref{sec:IGT}. Since $V_{\theta}$ is a smooth approximation of $V$, it makes sense to penalize deviations from~\eqref{eq:hjb} and~\eqref{eq:sensitivity} during training. Furthermore, given an initial guess $\mathrm{x}^0$ of an optimal trajectory for Problem~\eqref{def:oc_problem}, in view of~\eqref{eq:sensitivity} it is also meaningful to provide $[t,T]\ni s\mapsto (\mathrm{x}^0(s),\nabla_{x}V_{\theta}(s,\mathrm{x}^0(s)))\in\RR^{d}\times\RR^{d}$ as initial guess in the implementation of a numerical method to solve the TPBVP~\eqref{eq:TPBVP}.  
	
	\subsection{First-order MFG systems} \label{subsec:MFG}
	Consider now a population of agents distributed at time $t=0$ as $m_{0}\in\P(\RR^{d})$, where $\P(\RR^{d})$ denotes the space of probability measures over $\RR^{d}$. Given a curve of probability measures $[0,T]\ni t\mapsto m(t)\in\P(\RR^{d})$, where $m(t)$ represents the distribution of the agents at each time $t\in [0,T]$, a {\it typical agent} positioned at $x\in \RR^{d}$ at time $t\in [0,T]$ solves the optimal control problem~\eqref{def:oc_problem} with 
	\begin{equation}
		\ell(s,y,\alpha)=\ell_{0}(s,y,\alpha)+F(s,y,m(s))\quad\text{and}\quad g(y)=G(y,m(T))\quad\text{for all }(s,y,\alpha)\in [t,T]\times\RR^{d}\times A, 
		\label{eq:ell_mfg}
	\end{equation}
	where $\ell_{0}\colon[0,T]\times\RR^{d}\times A\to\RR$ is a cost not depending on $m$ and $F\colon[0,T]\times\RR^{d}\times\P(\RR^{d})\to\RR$ and $G\colon\RR^{d}\times\P(\RR^{d})\to\RR$ are known as the {\it coupling functions}. The associated value function $V[m]$ of a typical agent solves the HJB equation 
	\begin{equation}
		\begin{aligned}   
			-\partial_{t} V(t,x) +\widetilde{\H}(t,x,\nabla_{x}V(t,x))&= F(x,m(t))\quad\text{for }(t,x)\in ]0,T[\times \RR^d, \\
			V(T,x) &=G(x,m(T)) \quad\text{for }x\in  \RR^d,
		\end{aligned}
		\label{eq:hjb_m}
	\end{equation}
	where $\widetilde{\H}(t,x,p)=\sup_{\alpha\in A}\left\{-\ell_{0}(t,x,\alpha)-p\cdot b(t,x,\alpha)\right\}$ for every $(t,x,p)\in[0,T]\times\RR^{d}\times\RR^{d}$. Now, if $V[m]$ is sufficiently regular and agents located at $x$ at time $t=0$ provide their {\it best response} to $m$, i.e. act optimally satisfying~\eqref{eq:optimal_dynamics} with $t=0$, then the evolution of the distribution of the agents $[0,T]\ni t\mapsto {\rm BR}[m](t)\in\P(\RR^{d})$ is given by 
	\begin{equation}
		{\rm BR}[m](t)(\mathcal{O})=m_{0}\left(\left(\Phi[m](t,\cdot)\right)^{-1}(\mathcal{O})\right)\quad\text{for all }t\in [0,T],\,\mathcal{O}\in\mathcal{B}(\RR^{d}),
		\label{eq:push_forward_initial_measure}
	\end{equation}
	where, for every $(t,x)\in [0,T]\times\RR^{d}$, $\Phi[m](t,x)$ denotes the solution to 
	\begin{equation}
		\left\{\begin{aligned}
			\dot{\x}(s) &=-\partial_{p}\widetilde{\H}\left(s,\x(s),\nabla_{x}V[m](s,\x(s))\right)\quad\text{for a.e. }s\in [0,T],\\
			\x(0)&=x.
		\end{aligned}
		\right.
		\label{eq:optimal_dynamics_time_0}
	\end{equation}
	at time $t$. In other words, the {\it flow} $\Phi[m]$ satisfies 
	\begin{equation}
		\left\{
		\begin{aligned}
			\partial_{t}{\Phi[m]}(t,x) &=-\partial_{p}\widetilde{\H}\left(t,\Phi[m](t,x),\nabla_{x}V[m](t,\Phi[m](t,x))\right)\quad\text{for a.e. }t\in [0,T],\\
			\Phi[m](0,x)&=x. 
		\end{aligned}
		\right.
		\label{eq:Phi_EDP}
	\end{equation}
	Equivalently, standard results (see, e.g.,~\cite[Chapter 8]{ambrosio2008gradient} and the references therein) show that ${\rm BR}[m]$ is given by the solution to the following {\it continuity equation}:
	\begin{equation}
		\begin{aligned}
			\partial_{t}\mu-\text{div}\left(\partial_{p}\widetilde{\H}(t,x,\nabla_{x}V[m](t,x))\mu\right)&=0\quad\text{in }]0,T[\times \RR^d, \\
			\mu(0) &=m_{0}\quad\text{in }\P(\RR^{d}).
		\end{aligned}
		\label{eq:continuity_eq_depending_on_m}
	\end{equation}
	In this context, if $[0,T]\ni t\mapsto m^*(t)\in\P(\RR^{d})$ is such that $m^*={\rm BR}[m^*]$, one then says that $m^*$ is a MFG equilibrium. In turn, MFG equilibria are described by the following system of coupled equations:
	\begin{equation}
		\begin{aligned}   
			-\partial_{t} V(t,x) +\widetilde{\H}  (t,x, \nabla_{x}V(t,x))&= F(x,m(t))\quad\text{for }(t,x)\in ]0,T[\times \RR^d, \\
			V(T,x) &=G(x,m(T))\quad\text{for }x\in  \RR^d,\\
			\partial_{t}m-\text{div}\left(\partial_{p}\widetilde{\H}(t,x,\nabla_{x}V(t,x))m\right)&=0\quad\text{for }(t,x)\in ]0,T[\times \RR^d, \\
			m(0) &=m_{0}\quad\text{in }\P(\RR^{d}).
		\end{aligned}
		\label{eq:MFG_system}
	\end{equation}
	The existence of solutions to~\eqref{eq:MFG_system} has been studied  in~\cite{MR2295621,MR4214773,MR4304905,MR4102464} under several assumptions on the data $\ell_0$, $b$, and the coupling functions $F$ and $G$. The uniqueness of an equilibrium can be expected under structural assumptions on the data which often involve a monotonicity property of the cost in terms of the distribution of the agents. More precisely,  if the coupling functions  $\Xi(x,\mu)=F(x,\mu),\, G(x,\mu)$ satisfy
	\begin{equation}
		\ds \int_{\RR^d}\big(\Xi(x,\mu_1)-\Xi(x,\mu_2)\big)\dd(\mu_1-\mu_2)(x)  \geq 0 \quad \text{for all $\mu_1$, $\mu_2\in \P_1(\RR^d)$},
		\label{eq:monotonia_Phi}
	\end{equation}
	\smallskip 
	then one can show that the equilibrium $m^*$ is unique (see, e.g.,~\cite{MR2295621,These_Saeed_18,MR4835150}).  
	
	In some cases, the solution to~\eqref{eq:MFG_system} can be approximated through the so-called {\it fictitious play} method, introduced in the context of potential MFGs in~\cite{MR3608094} and later extended to first-order (non-potential) MFGs in~\cite{These_Saeed_18} (see also~\cite{MR4835150,MR4030259} for the application of this method to discretized versions of~\eqref{eq:MFG_system}).
	The iterates can be computed by solving separately the HJB and the continuity equations and read as follows:
	\begin{equation}
		\begin{aligned}
			\text{Given }\overline{m}_0\colon[0,T]\to\P(\RR^{d}),\\
			(\forall\,k\geq 0)\quad \mu_{k+1}={\rm BR}[\overline{m}_{k}], \;\; \overline{m}_{k+1}=\overline{m}_{k}+\frac{1}{k+1}(\mu_{k+1}-\overline{m}_{k}). 
		\end{aligned}
		\label{eq:fictitious_play_method}
	\end{equation}
	
	We refer the reader to~\cite[Chapter 3]{These_Saeed_18} for the study of the convergence of both sequences $(\mu_{k})_{k\in\NN}$ and $(\overline{m}_{k})_{k\in\NN}$ towards a MFG equilibrium in the space $C([0,T];\P_{1}(\RR^{d}))$, where $\P_{1}(\RR^{d})$ denotes the space of probability measures on $\RR^{d}$, with finite first-order moment, endowed with the $1$-Wasserstein distance (see, e.g.,~\cite[Chapter 7]{ambrosio2008gradient}). To monitor convergence, in our numerical tests in Section~\ref{subsec:numerical_mfg} we stop the iterates at $k$ as soon as the aggregated Sinkhorn divergence ${\bf S_{\varepsilon}^{\infty}}(\overline{m}_{k}(t),{\rm BR}[\overline{m}_{k}](t))$ (defined in Section \ref{sec:MFG_approximation} below) is smaller than a given tolerance. We also check that the so-called {\it exploitability} (see e.g.~\cite{Lauriere_fictitious_play}) associated with large iteration numbers is small. More precisely, given a smooth feedback control $\alpha\colon[0,T]\times\RR^{d}\to A$, its exploitability is defined as
	\begin{equation}
		\psi(\alpha)=\int_{\RR^{d}}\left(J[\mu^{\alpha}](x,\alpha)-V[\mu^{\alpha}](0,x)\right)\dd m_{0}(x),
		\label{eq:exploitability_continuous}
	\end{equation}
	where $\mu^{\alpha}$ is the solution to 
	\begin{equation}
		\begin{aligned}
			\partial_{t}\mu+\text{div}\left(b(t,x,\alpha(t,x))\mu\right)&=0\quad\text{in }]0,T[\times \RR^d, \\
			\mu(0) &=m_{0}\quad\text{in }\P(\RR^{d}),
		\end{aligned}
		\label{eq:continuity_eq_driven_by_alpha}
	\end{equation}
	and, for every $x\in\RR^{d}$ and $[0,T]\ni t\mapsto \mu(t)\in\P(\RR^{d})$,
	$$
	J[\mu](x,\alpha)=\int_{0}^{T}\ell(t,\x^{\alpha}(t),\alpha(t,\x^{\alpha}(t)),\mu(t))\dd t+G(\x^{\alpha}(T),\mu(T)),
	$$
	with $\ell$ given in~\eqref{eq:ell_mfg} and $\x^\alpha$ being the solution to
	\begin{equation}
		\label{eq:din_control_feedback}
		\left\{\begin{aligned}
			\dot{\x}(t) &= b\left(t,\x(t),\alpha(t,\x(t))\right)\quad\text{for a.e. }t\in [0,T],\\
			\x(0)&=x.
		\end{aligned}
		\right.
	\end{equation}
	Notice that $\psi(\alpha)\geq 0$ and $\psi(\alpha)=0$ if and only if $(V[\mu^{\alpha}],\mu^{\alpha})$ solves system~\eqref{eq:MFG_system}. Therefore, if convergence holds, one expects a small exploitability for large values of $k$.  
	
	\section{The initializing, generating, and training method to approximate optimal feedback controls}
	\label{sec:IGT}
	In this section we propose an approximation $V^{\mathrm{NN}}_{\theta}$ of the value function $V$, associated with the parameterized family of problems~\eqref{def:oc_problem}, through neural networks. In view of~\eqref{eq:optimal_feedbacks}, we obtain an approximation of optimal feedback controls taking the form 
	\begin{equation}
		[0,T]\times\RR^{d}\ni(t,x)\mapsto \Psi\big(t, x, \nabla_{x} V^{\mathrm{NN}}_{\theta}(t, x)\big)\in A.
		\label{def:approximated_feedback}
	\end{equation}
	The architecture of the neural networks used for this approximation is described below. As in~\cite{nakamura2021adaptive}, we train $V^{\mathrm{NN}}_{\theta}$ with a dataset obtained by solving in open-loop $({\bf P}^{t_{b},x_{b}})$ over a batch of points $\{(t_{b},x_{b})\}_{b=1}^{B}$. A key distinction from~\cite{nakamura2021adaptive} lies in the resolution of the associated TPBVP problems. Specifically, we construct an initial guess by using a rough approximation of 
	$V$ through the Deep Galerkin method~\cite{sirignano2018dgm} in order to foster global optimality. This strategy is inspired by~\cite{cristiani2010initialization}, where the authors approximate 
	$V$ by discretizing~\eqref{eq:hjb} with a finite difference scheme, which is feasible for low state dimensions, and employ the obtained approximation as an initial guess to solve~\eqref{eq:TPBVP} through the so-called shooting method.
	
	\vspace{0.2cm}
	
	{\it Neural network approximation}. We utilize multilayer feedforward NNs similar to those in previous methods (see, e.g.,~\cite{lin2021alternating}). Despite the availability of more complex architectures for other applications, we have designed a specific model adapted for efficient computation. We parameterize the value function as 
	\begin{equation}
		V^{\mathrm{NN}}_{\theta}(t, x) = \left( 1 - \varphi(t) \right) \mathrm{N}_{\theta}(t, x) + \varphi(t) g(x)\quad \text{for all } (t,x) \in [0,T] \times \RR^{d},
		\label{para hjb}
	\end{equation}
	where $\mathrm{N}_{\theta} \colon [0,T] \times \RR^{d} \rightarrow \RR$ is a neural network and $\varphi \colon [0,T] \to [0,1]$ is a smooth transition function satisfying $\varphi(T) = 1$. This construction guarantees that the terminal condition is exactly satisfied, i.e., $V^{\mathrm{NN}}_{\theta}(T, x) = g(x)$ for all $x \in \RR^{d}$. In practice, we consider two choices for the function $\varphi$:
	\begin{equation}
		\varphi_1(t) = \exp(t - T)\quad\text{and}\quad\varphi_2(t) = \frac{t}{T}\quad\text{for all }t \in [0,T]. 
		\label{eq:parametrizations_phi}
	\end{equation}
	Regarding $\mathrm{N}_{\theta}$, given $L,\,n\in\NN$, layers $h_{[1]}\colon [0,T]\times\RR^{d}\to\RR^{n}$, $h_{[i]}\colon\RR^{n}\to\RR^{n}$ ($i=2,\hdots,L-1$), and $h_{[L]}\colon\RR^{n}\to\RR$, we employ the following ResNet architecture:
	\begin{equation}
		\mathrm{N}_{\theta}(t, x)=\left(h_{[L]} \circ h_{[L-1]} \circ\ldots \circ h_{[1]}\right)(t,x)\quad\text{for all }(t,x)\in[0,T]\times\RR^{d}.
		\label{def:neural_network}
	\end{equation}
	Given smooth activation functions $\sigma_{i}\colon\RR\to\RR$ ($i=1,\hdots,L-1$) and connection weights $\beta_{i}\in]0,1[$ ($i=2,\hdots,L-1$), we consider layers of the form
	\begin{equation}
		\begin{aligned}
			h_{[1]}(z) &= \sigma_{1}(W_{1} z + c_{1})\quad\text{for all }z\in[0,T]\times\RR^{d},\\
			h_{[i]}(y)&=\sigma_{i}(y+\beta_{i}(W_{i}y + c_{i})) \quad \text{for all }i=2,\hdots,L-1,\, y\in\RR^{n},\\
			h_{[L]}(y)&=W_{L}y+c_{L}\quad\text{for all }y\in\RR^{n},
		\end{aligned}
    \label{eq:architecture_1}
	\end{equation}
where $W_{1}\in\RR^{n\times(1+d)}$, $W_{i}\in\RR^{n\times n}$ ($i=2,\hdots,L-1$), $W_{L}\in\RR^{1\times n}$, $c_{i}\in\RR^{n}$ ($i=1,\hdots,L-1$), $c_{L}\in\RR$, and the functions $\sigma_i$ are applied component-wise. In~\eqref{def:neural_network} the parameter $\theta=\{(W_{i},c_{i})\,|\,i=1,\hdots,L\}$ is the so-called trainable weight.  We denote by $\Theta$ the set of trainable weights.

\medskip
	
	{\it Initialization step.} To provide a suitable initial guess for solving the TPBVP~\eqref{eq:TPBVP}  through an iterative method, we first compute a rough estimate of the solution to the HJB equation~\eqref{eq:hjb}. This is achieved by employing the DGM~\cite{sirignano2018dgm}, which consists of training a neural network through the minimization of the loss function
	\begin{equation}
		\Theta \ni \theta \mapsto \lossHJB(\theta) := \frac{1}{M} \sum_{m=1}^{M} \left| \partial_t V^{\mathrm{NN}}_{\theta}(t_{m},x_{m})
		- \H\left(t_{m},x_{m}, \nabla_{x}V^{\mathrm{NN}}_{\theta}(t_{m},x_{m})\right) \right|^2 \in \RR,
		\label{def:lossHJB_init}
	\end{equation}
	thus providing a first approximation $V^{\mathrm{NN}}_{\theta_0}$. Here, $\{(t_{m}, x_m)\}_{m=1}^{M}$ denotes a set of randomly sampled points in the domain $[0,T] \times \mathbb{R}^d$.  In complex cases where solving the HJB equation over the entire domain with DGM becomes difficult, we adopt a characteristic-driven DGM approach (C-DGM). This method can be seen as a simplified version of the ML method proposed in~\cite{onken2022neural},  which will be sufficient in practice to initialize the data generation process in the next step. Specifically, for $\theta\in\Theta$ define the feedback control
	\begin{equation}
		\alpha_{\theta}(t,x)=\Psi(t,x,\nabla_{x}V^{\mathrm{NN}}_{\theta}(t,x))\quad\text{for all }(t,x)\in[0,T]\times\RR^{d}. 
		\label{eq:feedback_control_NN}
	\end{equation}
	Given a first guess $\theta_0\in\Theta$ and an initial point $x_0\in\RR^{d}$,  we solve the ODE
	\begin{equation}
		\label{eq:din-avec-V-NN}
		\left\{\begin{aligned}
			\dot{\x}(t) &= b\left(t,\x(t),\alpha_{\theta_0}(t,\x(t))\right)\quad\text{for a.e. }t\in [0,T],\\
			\x(0)&=x_0
		\end{aligned}
		\right.
	\end{equation}
	by using  an Initial Value Problem (IVP) solver. In our implementations, we use Python's \texttt{solve\_ivp} function, which yields a sequence of time points $t_1,\hdots,t_{N}\in (0,T]$ along with their corresponding approximate states $x_1,\hdots,x_{N}$, ensuring a specified level of precision. Next,  the HJB residual loss~\eqref{def:lossHJB_init} is computed over the collection of points $\{(t_{n}, x_n)\}_{n=1}^{N}$ and then $\theta_0$ is updated.  Additional updates of $\theta_0$ can be obtained by repeating this procedure for several initial conditions.

	\medskip
	
	{\it Data generation}.  To solve in open-loop the optimal control problem over a batch of initial times and initial states, we first solve $({\bf P}^{0,x_{0}^{i}})$ over a sample of initial states $\{x_{0}^{i}\}_{i=1}^{S}\subset\RR^{d}$. As explained below, the resolution of the TPBVP system then provides a family of time points and corresponding states, which are then used to train the neural network. This approach allows the network to be trained on dynamically generated data points, ensuring better coverage of the solution space. 
	
	Notice that if $(\x^*,\alpha^*)$ solves~$({\bf P}^{0,x_{0}^{i}})$, then $v(t):=V(t,\x^*(t))=\int_{t}^{T}\ell(s,\x^{*}(s),\alpha^*(s))\dd s+g(\x^*(T))$ satisfies 
	$$
	\left\{
	\begin{aligned}
		\dot{v}(t)&=-\ell(t,\x^{*}(t),\alpha^*(t))\quad\text{for a.e. }t\in [0,T],\\
		v(T) &=g(\x^*(T)).
	\end{aligned}
	\right.
	$$
	For numerical purposes, it will be useful to add this equation to the TPBVP  system~\eqref{eq:TPBVP} which, recalling~\eqref{eq:unique_minimizer_hamiltonian}, yields
	\begin{equation}
		\left\{\begin{aligned}
			\dot{v}(t)&=-\ell(t,\x^{*}(t),\Psi(t,\x^*(t),p^{*}(t)))\quad\text{for a.e. }t\in [0,T],\\
			\dot{\x}^*(t) &= b\left(t, \x^*(t), \Psi(t,\x^*(t),p^{*}(t))\right)\quad\text{for a.e. }t\in [0,T],\\
			\dot{p}^*(t) &= \nabla_{x}H\left(t, \x^*(t), p^*(t), \Psi(t,\x^*(t),p^{*}(t))\right)\quad\text{for a.e. }t\in [0,T],\\
			0&= B((v(0),\x^*(0),p^*(0)),(v(T),\x^*(T),p^*(T)))
		\end{aligned}
		\right.
		\label{eq:TPBVP_augmented}  
	\end{equation}
	where $B\colon(\RR\times\RR^{d}\times\RR^{d})^{2}\to\RR\times\RR^{d}\times\RR^{d}$ is the function imposing the boundary conditions, defined as 
	$$
	B((v_0,x_0,p_0), (v_T,x_T,p_T))= (v_{T}-g(x_T),x_0-x_{0}^{i},p_T-\nabla g(x_T))
	$$
	for all $(v_0,x_0,p_0),\, (v_T,x_T,p_T)\in\RR\times\RR^{d}\times\RR^{d}$. System~\eqref{eq:TPBVP_augmented} corresponds to the optimality condition of problem~$({\bf P}^{0,x_{0}^{i}})$ and solving it will generate data comprising value function evaluations, state and adjoint state trajectories. 
	
In our numerical experiments, system~\eqref{eq:TPBVP_augmented} is solved by using Python's \texttt{solve\_bvp} function, which is a robust solver for TPBVP, based on the work~\cite{kierzenka2001bvp}, and has been applied to optimal control problems in~\cite{kang2017mitigating,nakamura2021adaptive}. It employs a 4th-order collocation method with residual control, using either piecewise cubic polynomials or an implicit Runge-Kutta formula, and solves the resulting finite dimensional system by employing Newton's method. In turn, as typical for TPBVP solvers, the method is highly sensitive to the choice of the initial guess,  especially when the interval $[0,T]$ is large. The rough estimate $V_{\theta_0}^{\rm NN}$ of the value function obtained in the initialization step is used to construct a suitable initial guess. More precisely, for every $i=1,\hdots,S$, we solve the ODE~\eqref{eq:din-avec-V-NN} with $x_0$ replaced by $x_{0}^{i}$  by using the IVP solver \texttt{solve\_ivp} to obtain 
	a sequence of time points $\tilde{t}_1^{i}, \hdots,\tilde{t}_{N^i}^{i}\in (0,T]$ along with their corresponding approximate states $\tilde{x}_1^{i},\hdots,\tilde{x}_{N^i}^{i}$. Setting $\tilde{t}_0^{i}=0$, $\tilde{x}_0^{i}=x_0^{i}$ and defining, for every $k=0,\hdots,N^i$,  $\tilde{v}_{k}^{i}=V^{\mathrm{NN}}_{\theta_0}(\tilde{t}_k^{i},\tilde{x}_k^i)$ and $\tilde{p}_{k}^{i}=\nabla_{x}V^{\mathrm{NN}}_{\theta_0}(\tilde{t}_k^{i},\tilde{x}_k^i)$, we use $\{(\tilde{v}_{k}^{i},\tilde{x}_{k}^{i},\tilde{p}_{k}^{i})\}_{k=0}^{N_i}$ as initial guess to solve the TPBVP~\eqref{eq:TPBVP_augmented}, yielding at last a dataset
	\begin{equation}
		\mathcal{D}_{\text{OC}}=\left\{\left(t^{b}, x^{b}, v^{b}, p^{b}\right)\right\}_{b=1}^{B}
		\label{eq:data_ini}
	\end{equation}
	which will play an essential role for training the model in the following step.
	
	{\it Model training.} Having at our disposal the data set $\mathcal{D}_{\text{OC}}$, we train the NN by minimizing a suitable loss function. More precisely, to update $\theta$ we consider the following problem:
	\begin{equation} 
		\min_{\theta\in\Theta}\;\mathrm{Loss}_{\mathrm{val}}(\theta),
		\label{eq:minimization_complete_loss}
	\end{equation}
	where 
	\begin{equation}
		\mathrm{Loss}_{\mathrm{val}}(\theta):=\lossv(\theta)+\lossgradv(\theta)+ \lambda_1\lossHJB(\theta)\quad\text{for all }\theta\in\Theta,
		\label{eq:loss_total}
	\end{equation}
	with $\lambda_{1}>0$, $\lossHJB$ being defined by~\eqref{def:lossHJB_init}, and
	\begin{equation}\label{eq:loss_V_grad_HJB}
		\begin{aligned}
			&\Theta\ni\theta\mapsto\lossv(\theta):=\frac{1}{B} \sum_{b=1}^{B}\left|v^{b}-V_{\theta}^{\mathrm{NN}}\left(t^{b}, x^{b} \right)\right|^2\in\RR,
			\\
			&\Theta\ni\theta\mapsto\lossgradv(\theta):=\frac{1}{B} \sum_{b=1}^{B}\left\|p^{b}-\nabla_{x} V^{\mathrm{NN}}_{\theta}\left(t^{b}, x^{b} \right)\right\|^2\in\RR,
		\end{aligned}
	\end{equation}
	penalize, respectively, deviations from the estimates of the value function and its gradient. Depending on the application at hand, the term $\lossHJB(\theta)$ can be computed over the batch $\{(t^b,x^b)\}_{b=1}^{B}$ from $\mathcal{D}_{\text{OC}}$ or over a different batch of time-space points to enforce global optimality. We will frequently adopt this last variant to deal with problems where the state dimension is not very large. 
	
	In our simulations, Problem~\eqref{eq:minimization_complete_loss} is solved by using the stochastic gradient method, and provides an update of the value of $\theta$, which can then be used in a second round to initialize the solution to TPBVP~\eqref{eq:TPBVP_augmented} with a better initial guess. The procedure can thus be performed iteratively for a number $R$ of rounds to produce a final and sharp  estimation of the value function, of its gradient, and, in turn, of the approximate optimal feedback control~\eqref{def:approximated_feedback}. In our numerical experiments rather sharp results are already obtained for small values of $R$ (typically $R\leq 5$). The method, which we call Initializing, Generating, and Training (IGT), is summarized in Algorithm~\ref{alg:IGT} below.
	
	\begin{algorithm}
		\begin{algorithmic}
			\Require  Batch sizes $M$ and $S$, number of rounds $R$.
			\Require  Initialize neural network parameters $\theta\in\Theta$.
			\smallskip
			\State {\it Initialization} 
			\smallskip
			\While{not converged\footnotemark}
			\State (if DGM) Sample a batch $\{(t_{m},x_{m})\}_{m=1}^{M}$ from  $[0, T] \times \mathbb{R}^d$.
			\State (if C-DGM) Generate  $\{(t_{n},x_{n})\}_{n=1}^{N}$ by solving~\eqref{eq:din-avec-V-NN}.
			\State Compute $\lossHJB$ with~\eqref{def:lossHJB_init}.
			\State Backpropagate $\lossHJB$ to update $\theta$.
			\EndWhile
			\smallskip
			\For{$r = 1$ \textbf{to} $R$}
			\State {\it Data generation}
			\State Generate $\mathcal{D}_{\text{OC}}$, given by~\eqref{eq:data_ini}, using a batch of initial conditions $\{x_{0}^{i}\}_{i=1}^{S}$, $V^{\mathrm{NN}}_{\theta}$, and $\nabla_{x}V^{\mathrm{NN}}_{\theta}$.
			\smallskip    
			\State {\it Training}
			\While{not converged}
			\State Compute $\mathrm{Loss}_{\mathrm{val}}$, defined in~\eqref{eq:loss_total}, with data $\mathcal{D}_{\text{OC}}$.
			\State Backpropagate $\mathrm{Loss}_{\mathrm{val}}$ to update $\theta$.
			\EndWhile
			\EndFor
			\Statex
			\Return $\theta$
		\end{algorithmic}
		\caption{IGT Algorithm}
		\label{alg:IGT}
	\end{algorithm}
	\footnotetext{ In this context, {\it converged} refers to $V_{\theta}^{\mathrm{NN}}$ providing an initialization for which the solver of the TPBVP in the data generation step converges.
	}

		\section{Approximation of solutions to high-dimensional first-order mean field games systems}
		\label{sec:MFG_approximation}
		The aim of this section is to approximate the equilibria of the MFG problem introduced in Subsection~\ref{subsec:MFG}. As is standard in the MFG literature, we rely on the fictitious play method~\eqref{eq:fictitious_play_method} to construct an approximation of the equilibrium. In this framework, the choice of the initial distribution plays a crucial role, since the iterates are defined through successive averaging of past distributions. For this reason, we introduce an initialization procedure designed to provide a more accurate starting distribution for the fictitious play algorithm, thereby improving the quality and stability of the resulting approximation.
		
		Throughout this section, we consider two neural network approximations, denoted by $V^{\mathrm{NN}}_{\theta}$ and $\Phi^{\mathrm{NN}}_{\omega}$, which approximate, respectively, the value function of the optimal control problem associated with the MFG system and the flow map inducing the best-response operator.
		
		More precisely, the proposed method to approximate MFG equilibria proceeds in two stages. In the first stage, we seek $\overline{m}_0\colon [0,T]\to \mathcal{P}(\RR^d)$ to initialize the fictitious play iterations \eqref{eq:fictitious_play_method}. This is achieved by jointly training two neural networks, namely the value network $V^{\mathrm{NN}}_{\theta_0}$ and the generator network $\Phi^{\mathrm{NN}}_{\omega_0}$.
		This initialization is carried out by minimizing a composite loss function that penalizes both the residual of the HJB equation~\eqref{eq:hjb_m} and the discrepancy between the generator dynamics and the optimal flow~\eqref{eq:Phi_EDP}. Then given $\{x_0^m\}_{m=1}^{M'}$ sampled from $m_0$ we define $\overline{m}_0(t)
		= \frac{1}{M'}\sum_{m=1}^{M'}\delta_{\Phi_{\omega_0}^{\text{NN}}(t,x_{0}^{m})}\in\mathcal{P}(\RR^{d})$ for $t\in [0,T]$.
		The precise formulation of this initialization step is given below.
		
		In the second stage, starting from the initial distribution $\overline{m}_0$, the approximation is iteratively refined using the fictitious play method~\eqref{eq:fictitious_play_method}. At each iteration, for a given mean field flow $\overline{m}$, the IGT method presented in the previous section is employed to efficiently approximate the solution $V[\overline{m}]$ of~\eqref{eq:hjb_m}. The resulting approximation of the feedback control~\eqref{eq:feedback_control_NN} is then used to construct a neural network approximation of the flow map $\Phi[\overline{m}]$ defined by~\eqref{eq:Phi_EDP}. This yields an approximation of the best response ${\rm BR}[\overline{m}]$ as the empirical measure induced by the learned flow map. 
		Combining both approximations with the fictitious play iterates, we obtain a ML method to approximate solutions to high-dimensional first-order MFG systems. Although the IVP solver already produces the trajectories used for training, the resulting data is discrete. The generator instead approximates the continuous map $(t,x_0)\mapsto\Phi(t,x_0)$. Once trained, it can be evaluated inexpensively at arbitrary times and for new samples from $m_0$, without solving a new IVP for every evaluation. This  representation is useful when the same flow must be queried repeatedly to approximate coupling terms, empirical distributions, and convergence indicators. For a fixed fictitious-play iteration, the generator also replaces a potentially dense trajectory data set by a fixed collection of network parameters.

	{\it  Neural network approximation.} The mapping $\Phi^{\mathrm{NN}}_{\omega}$ is implemented as a multilayer feedforward neural network, analogous to $V^{\mathrm{NN}}_{\theta}$, with outputs taking values in $\mathbb{R}^d$, where $V^{\mathrm{NN}}_{\theta}$ follows the same architecture as in the previous section. 
	Here, $\omega$ denotes a parameter belonging to a parameter space $\mathcal{W}$ and the architecture used to build each coordinate of $\Phi_{\omega}^{\text{NN}}$ is similar to the one used for $V_{\theta}^{\text{NN}}$, employing a smooth transition function to ensure the equality $\Phi_{\omega}^{\text{NN}}(0,x)=x$ for all $x\in\RR^{d}$. 
		
	{\it  Initialisation.} To initialize the iterative procedure, we jointly train a value network $V_{\theta}^{\text{NN}}$ and a generator network $\Phi_{\omega}^{\text{NN}}$ in order to compute $\overline{m}_0$ as described above. 
	Both networks are trained simultaneously by minimizing the loss function
	\begin{equation}
		\mathcal{L}(\theta, \omega) = \lossHJB^0(\theta, \omega) +  \lossODE^0(\theta, \omega),
		\label{eq:joint_loss_mfg}
	\end{equation}

	The term $\lossHJB^0$ penalizes the residual of the HJB equation~\eqref{eq:hjb_m} and is defined as
	\begin{equation}
		\begin{aligned}
			\lossHJB^0(\theta, \omega) := & \frac{1}{M_1'} \sum_{m=1}^{M_1'} \left| \partial_t V_{\theta}^{\text{NN}}(t_m, x_m) - \widetilde{\H}\left(t_m, x_m, \nabla_x V_{\theta}^{\text{NN}}(t_m, x_m)\right) + F(x_m, \overline{\mu}_{\omega}(t_m)) \right|^2 \\
			& + \frac{1}{M_2'} \sum_{m=1}^{M_2'} \left| \partial_t V_{\theta}^{\text{NN}}(\tau_m, y_m) - \widetilde{\H}\left(\tau_m, y_m, \nabla_x V_{\theta}^{\text{NN}}(\tau_m, y_m)\right) + F(y_m, \overline{\mu}_{\omega}(\tau_m)) \right|^2,
		\end{aligned}
		\label{eq:loss_hjb_mfg_init}
	\end{equation}
	where $\{(t_m, x_m)\}_{m=1}^{M_1'}$ are a set of randomly sampled points in the domain $[0,T] \times \mathbb{R}^d$, and $y_m = \Phi_{\omega}^{\text{NN}}(\tau_m, z_0^m)$ are points along the trajectories generated by the flow map, with $\{\tau_m\}_{m=1}^{M_2'}$ and $\{z_{0}^{m}\}_{m=1}^{M_2'}$ independent samples from the uniform distribution in $[0,T]$ and $m_0$, respectively. Here $\overline{\mu}_{\omega}(t) := \frac{1}{M'}\sum_{m=1}^{M'}\delta_{\Phi_{\omega}^{\text{NN}}(t,x_{0}^{m})}\in\mathcal{P}(\RR^{d})$ denotes the empirical measure associated with the flow $\Phi_{\omega}^{\text{NN}}$, with $\{x_0^m\}_{m=1}^{M'}$ sampled from $m_0$.
	The term $\lossODE^0$ enforces the consistency of the generator network with the optimal dynamics~\eqref{eq:Phi_EDP} and is given by
	\begin{equation}
		\lossODE^0(\theta, \omega) := \frac{1}{M_3'} \sum_{m=1}^{M_3'} \left\| \partial_t \Phi_{\omega}^{\text{NN}}(s_m, \xi_0^m) - b\left(s_m, \Phi_{\omega}^{\text{NN}}(s_m, \xi_0^m), \alpha_{\theta}(s_m, \Phi_{\omega}^{\text{NN}}(s_m, \xi_0^m))\right) \right\|^2,
		\label{eq:loss_ode_mfg_init}
	\end{equation}
	where $\{s_{m}\}_{m=1}^{M_3'}$ and $\{\xi_{0}^{m}\}_{m=1}^{M_3'}$ are independent samples from the uniform distribution in $[0,T]$ and $m_0$, respectively, and $\alpha_{\theta}$ is the feedback control associated with $V_{\theta}^{\text{NN}}$ (see~\eqref{eq:feedback_control_NN}).
		This pre-training step yields two neural networks, $V_{\theta_0}^{\text{NN}}$ and $\Phi_{\omega_{0}}^{\text{NN}}$, the latter is used to define the initial distribution defined as $\overline{m}_0(t) = \overline{\mu}_{\omega_0}(t)$ for all $t\in [0,T]$.

		{\it  Data generation.} Given $\overline{m}$,  let $V_{\theta^*}^{\text{NN}}[\overline{m}]$ be the approximation of $V[\overline{m}]$ obtained by the IGT method.  In order to train the generator network we first consider a sample of initial states $\{x_{0}^{i}\}_{i=1}^{S}$ drawn from $m_0$ and, for every $i=1,\hdots,S$, we use an IVP solver to generate an approximation of the solution to 
		\begin{equation}
			\left\{\begin{aligned}
				\dot{\x}(s) &=b(s,\x(s),\alpha_{\theta^*}(s,\x(s)))\\
				&=-\partial_{p}\widetilde{\H}\left(s,\x(s),\nabla_{x}V_{\theta^*}^{\text{NN}}[\overline{m}](s,\x(s))\right)\quad\text{for a.e. }s\in [0,T],\\
				\x(0)&=x_{0}^{i},
			\end{aligned}
			\right.
			\label{eq:optimal_dynamics_time_0_with_VNN}
		\end{equation}
		where we recall that $\alpha_{\theta^*}$ is defined by~\eqref{eq:feedback_control_NN}. This yields  a sequence of times $t_{1}^{i},\hdots,t_{N^{i}}^{i}$ and states $x_{1}^{i},\hdots,x_{N^{i}}^{i}$, ensuring a required level of precision, producing the data set 
		\begin{equation}
			\mathcal{D}_{\text{MFG}}=\left\{\left(t_{j}^{i},x_{j}^{i}\right)\,|\,i=0,\hdots,S, \, j=1,\dots,N^{i}\right\}.
			\label{eq:data_G}
		\end{equation} 
		
		{\it  Model training.} The neural network $\Phi_{\omega}^{\text{NN}}[\overline{m}]$ is trained by penalizing deviations from the data set $\mathcal{D}_{\text{MFG}}$ and the residual of the ODE dynamics~\eqref{eq:Phi_EDP}, with $V[\overline{m}]$ replaced by $V^{\text{NN}}_{\theta^*}[\overline{m}]$, over points $\{(t^m,x_{0}^{m})\}_{m=1}^M$, where $\{t_{m}\}_{m=1}^{M}$ and $\{x_{0}^{m}\}_{m=1}^{M}$ are independent samples from the uniform distribution in $[0,T]$ and $m_0$, respectively. More precisely, to update $\omega$ one considers the following optimization problem \begin{equation} 
			\min_{\omega\in\mathcal{W}}\;\mathrm{Loss}_{\mathrm{gen}}(\omega),
			\label{eq:minimization_complete_loss_generator}
		\end{equation} 
		where 
		\begin{equation}
			\mathrm{Loss}_{\mathrm{gen}}(\omega):= \lossg(\omega)+\lambda_{2}\lossODE(\omega)\quad\text{for all }\omega\in\mathcal{W},
			\label{eq:loss_total_gen}
		\end{equation}
		with $\lambda_{2}>0$ and, for every $\omega\in\mathcal{W}$,
		\begin{equation}
			\begin{aligned}
				&\lossg(\omega):=\frac{1}{S} \sum_{i=1}^S\frac{1}{N_{i}}\sum_{j=1}^{N_{i}}\left\|x_j^i-\Phi_{\omega}^{\mathrm{NN}}[\overline{m}]\left(t_j^i,x_0^i \right)\right\|^2,\\
				&\lossODE(\omega):=\frac{1}{M} \sum_{m=1}^{M}\Big\|\partial_{t} \Phi_{\omega}^{\mathrm{NN}}[\overline{m}](t_{m},x_{0}^{m})- b\left(t_m,\Phi^{\mathrm{NN}}_{\omega}[\overline{m}](t_{m},x_{0}^{m}),\alpha_{\theta^*}(t_m, \Phi_{\omega}^{\mathrm{NN}}[\overline{m}](t_{m},x_{0}^{m}))\right)\Big\|^2.
			\end{aligned}
			\label{eq:losses_generator}
		\end{equation}
		The last update of $\omega$ after convergence is achieved is denoted as $\omega^*$. 
		\medskip
		
		{\it  Update and error measures.} Given the parameter $\omega^*\in\mathcal{W}$ obtained in the previous step, we consider a new batch $\B=\{x_{0}^{b}\}_{b=1}^{B}$ of initial conditions sampled from $m_{0}$ and approximate $\text{BR}[\overline{m}]$ by the following curve of empirical measures (see \eqref{eq:push_forward_initial_measure})
		\begin{equation}
			[0,T]\ni t\mapsto \overline{\mu}(t):=\frac{1}{B}\sum_{b=1}^{B}\delta_{\Phi_{\omega^*}^{\text{NN}}[\overline{m}](t,x_{0}^{b})}\in\mathcal{P}(\RR^{d}).
			\label{eq:mu_FP}
		\end{equation}
		Assuming that $\overline{m}$ has been already updated $k$ times, the $k+1$ update is defined through the fictitious play iteration~\eqref{eq:fictitious_play_method}:
		\begin{equation}
			(\forall\,t\in[0,T])\quad \overline{m}(t)\leftarrow \overline{m}(t)+\frac{1}{k+1}(\overline{\mu}(t)-\overline{m}(t)).
			\label{eq:fictitious_play_method_numerical}
		\end{equation}
		In our numerical experiments, we adopt a stopping criterion based on the proximity between the distributions $\overline{m}$ and $\overline{\mu}$ evaluated over a discrete time grid $\T \subset [0,T]$. To quantify this proximity, at each time $t\in\T$ we employ the Sinkhorn divergence $S_{\varepsilon}(\overline{m}(t),\overline{\mu}(t))$, as introduced in~\cite{genevay2018learning}, where $\varepsilon > 0$ denotes a regularization parameter. This divergence offers a computational advantage over classical optimal transport, as it can be evaluated more efficiently and exhibits significantly lower computational complexity, particularly in high-dimensional cases~\cite{genevay2019sample}. Specifically, we consider the following aggregated metric:
		\begin{equation}
			{\bf S_{\varepsilon}^{\infty}}(\overline{m},\overline{\mu}) := \max_{t\in\T} S_{\varepsilon}(\overline{m}(t), \overline{\mu}(t)),
			\label{eq:w_metrics}
		\end{equation}
		which is computed using the \texttt{geomloss} library~\cite{feydy2019interpolating}.  
		
		We also monitor the convergence to $0$ of the following approximation of the exploitability $\psi(\alpha_{\theta^*})$ (see~\eqref{eq:exploitability_continuous}):
		\begin{equation}
			\psi^{\B}(\alpha_{\theta^*})=\frac{1}{B}\sum_{b=1}^{B}\left(\tilde{J}[\overline{\mu}](x_{0}^{b},\alpha_{\theta^*})-V^{\text{NN}}_{\theta^*}[\overline{\mu}](0,x_{0}^{b})\right),
			\label{eq:exploitability_numerical}
		\end{equation}
		where, for a given uniform time grid $\{t_{k}\}_{k=0}^{N}\subset[0,T]$ with time step $\Delta t$, 
		$$
		\tilde{J}[\overline{\mu}](x_{0}^{b},\alpha_{\theta^*}):=\Delta t\sum_{k=0}^{N-1}\ell\left(t_k,\Phi_{\omega^*}^{\mathrm{NN}}[\overline{m}](t_{k},x_{0}^{b}),\alpha_{\theta^*}(t_k,\Phi_{\omega^*}^{\mathrm{NN}}[\overline{m}](t_{k},x_{0}^{b})),\overline{\mu}(t_k)\right)+G\left(\Phi_{\omega^*}^{\mathrm{NN}}[\overline{m}](T,x_{0}^{b}),\overline{\mu}(T)\right).
		$$
		
		From equation \eqref{eq:fictitious_play_method_numerical}, it follows that computing $\overline{m}$ at iteration $k$ of the fictitious play method requires using all the generators computed up to that iteration.
		Since storing and managing this growing collection of generators becomes increasingly expensive, to keep the method computationally feasible, we impose a maximum number of iterations $K_{{\rm max}}$ in the fictitious play procedure.
		We summarize in Algorithm~\ref{alg: IGT-MFG} the proposed procedure to approximate MFG equilibria. 
		
		\begin{algorithm}[hbt!]
			\begin{algorithmic}
				\Require Batch sizes $M_1',  \, M_2', \,M_3', \,M_{1},\,M_{2},\,S_{1},\,S_{2}$, number of rounds $R$, and maximum number of fictitious play iterations $K_{\rm{max}}$.
				\Require Batch of initial conditions $\B=\{x_{0}^{b}\}_{b=1}^{B}$ sampled from $m_0$ and  tolerance parameter $\texttt{tol}>0$. 
				\Require Initialize neural network parameters $\theta^{*}\in\Theta$ and $\omega^{*}\in\W$. 
				\Ensure $\delta_0\gets\texttt{tol}+1$.
				\State {\it Initialization}
				\While{not converged}
				\State Sample batches for $\lossHJB^0$ and $\lossODE^0$ computation. 
				\State Compute $\mathcal{L}(\theta^*, \omega^*)$ via~\eqref{eq:joint_loss_mfg}.
				\State Backpropagate $\mathcal{L}(\theta^*, \omega^*)$ to update $\theta^*$ and $\omega^*$.
				\EndWhile
				\State Set $\overline{m}_0(t) \leftarrow \overline{\mu}_{\omega^*}(t)$ for all $t \in [0,T]$.
				\smallskip
				\For{$k = 0$ \textbf{to} $K_{\max}-1$} 
				\If{$\delta_{k}\le\mathtt{tol}$} \textbf{break} \EndIf
				\State Update $\theta^*$ by approximating $V[\overline{m}_{k}]$ using the IGT method with input parameters $M_{1},\,S_1,\,R$.  
				\State {\it Data generation}
				\State \quad \quad  Generate $\mathcal{D}_{\text{MFG}}$, given by~\eqref{eq:data_G},  using a batch $\{x_{0}^{i}\}_{i=1}^{S_{2}}$, sampled from $m_0$, and $\nabla_{x}V^{\mathrm{NN}}_{\theta^{*}}[\overline{m}^k]$.
				\State {\it Training of $\Phi^{\rm NN}_{\omega}[\overline{m}_k]$}
				\While{not converged} 
				\State \quad \quad Sample a batch $\{x_{0}^{m}\}_{m=1}^{M_{2}}$ from  $m_{0}$ and $\{t_{m}\}_{m=1}^{M_{2}}$ from a uniform distribution in $[0,T]$.
				\State \quad \quad Compute $\text{Loss}_{\text{gen}}$, using $\mathcal{D}_{\text{MFG}}$ for $\text{Loss}_{\Phi}$ and the batch $\{(t_{m},x_{0}^{m})\}_{m=1}^{M_{2}}$ for $\lossODE$. 
				\State \quad \quad Backpropagate $\text{Loss}_{\text{gen}}$ to update $\omega^{*}$.
				\EndWhile
				\State {\it Update}
				\State \quad \quad  Compute $\mu_{k+1}$ using \eqref{eq:mu_FP} with the batch $\B$ and $\Phi^{\rm NN}_{\omega^*}[\overline{m}_{k}]$.
				\State \quad \quad  Update $\delta_{k+1} \gets {\bf S_{\epsilon}^{\infty}}(\overline{m}_{k},\mu_{k+1})$.
				\State \quad \quad  Update $\overline{m}_{k+1}\gets\overline{m}_{k}+\frac{1}{k+1}(\mu_{k+1}-\overline{m}_{k})$.
				\EndFor
				\State \Return $\theta^{*}$ and $\omega^{*}$.
			\end{algorithmic}
			\caption{IGT-MFG}
			\label{alg: IGT-MFG}
		\end{algorithm}

		\section{Numerical Results}
		\label{sec:numerical_results}
		In this section, we evaluate the performance of the proposed IGT and MFG-IGT algorithms. The methods are applied to a variety of optimal control and MFG problems to demonstrate their efficiency and robustness. 
		
			{\it Common implementation details.}
			Unless otherwise specified in an individual experiment, the value function network uses three hidden layers of width $128$, the residual-connection weight is $\beta_i=0.5$, and the activation function is $\tanh$. The generator network uses the analogous architecture with ReLU activations. Experiment-specific batch sizes, loss weights, and numbers of IGT rounds are reported with the corresponding numerical tests below.
		
		The implementation of all experiments is publicly available at~\url{https://github.com/Mouhcine-56/IGT}.
		In the tests admitting explicit expressions for the exact value function $V$, we evaluate the performance of the methods in terms of their relative $L^\infty$ and $L^2$ errors, defined, at time $t\in[0,T]$, as 
		$$
		E_{\infty}(V^{\mathrm{NN}}_{\theta}(t,\cdot)):= \frac{\max_i |V^{\mathrm{NN}}_{\theta}(t,x_i) - V(t,x_i)|}{\max_i |V(t,x_i)|}\quad\text{and}\quad E_{2}(V^{\mathrm{NN}}_{\theta}(t,\cdot)):=\frac{\left( \sum_{i} |V^{\mathrm{NN}}_{\theta}(t,x_i) - V(t,x_i)|^2 \right)^{1/2}}{\left( \sum_{i} |V(t,x_i)|^2 \right)^{1/2}}.
		$$
		
		To compare  optimal feedback controls $\alpha(t,x)\in\mathbb{R}^d$ with  the approximations $\alpha_\theta(t,x)\in\mathbb{R}^d$ defined in \eqref{eq:feedback_control_NN}, we employ the standard Euclidean norm. For each sampled point $x_i$, the pointwise control error is calculated as 
		$$
		\| \alpha_{\theta}(t,x_i) - \alpha(t,x_i)\|_2
		= \left( \sum_{j=1}^d 
		|\alpha_{\theta,j}(t,x_i) - \alpha_j(t,x_i)|^2 
		\right)^{1/2}.
		$$
		The associated relative $L^\infty$ and $L^2$ errors over the batch are defined as
		$$
		E_{\infty}(\alpha_{\theta}(t,\cdot))
		:= 
		\frac{\max_i \|\alpha_{\theta}(t,x_i) - \alpha(t,x_i)\|_2}
		{\max_i \|\alpha(t,x_i)\|_2},
		\qquad
		E_{2}(\alpha_{\theta}(t,\cdot))
		:=
		\frac{\left( \sum_i \|\alpha_{\theta}(t,x_i) - \alpha(t,x_i)\|_2^2 \right)^{1/2}}
		{\left( \sum_i \|\alpha(t,x_i)\|_2^2 \right)^{1/2}}.
		$$
		Here, the maximum and the sum are taken over a batch of sampled points $x_i$ which is made explicit for each test. The metrics above are computed at various choices of $t$ to assess both global and local accuracy of the neural network approximations.
		\subsection{Optimal control problems}
		
		In this section, we implement the IGT algorithm on a series of benchmark problems. We begin by studying a linear-quadratic  problem and a problem featuring a non-smooth value function with a known analytical solution. These examples enable a performance comparison between our method, the Adaptive Sampling and Model Refinement (ASMR) algorithm~\cite{nakamura2021adaptive}, and the Neural Network Approach (NNA)~\cite{onken2022neural}. We then turn to the quadcopter problem, a well-known example with practical relevance in real-world applications. Finally, we analyze an obstacle problem. For both the first and the last examples, we consider settings with dimensions $2$, $10$, $50$ and $100$ to evaluate the scalability of our approach.
		
		\subsubsection{{\bf Evaluating IGT}}\label{Analytic Comparison} {\it Linear-quadratic  problem.}
		To evaluate the effectiveness of the IGT algorithm, we consider first a very simple linear-quadratic example with an explicit solution which serves as a benchmark for numerical analysis. We have chosen to systematically compare the performance and reliability of the IGT, ASMR, and NNA methods as they share similar designs. These comparisons highlight the importance of the  initialization step, the penalization of the residual of the HJB equation, and the implementation of several rounds, each constituting a distinctive feature of the proposed IGT method.  
		
		We consider the family of optimal control problems~\eqref{def:oc_problem} with $T=1$, $A=\RR^{d}$, and 
		$$
		\ell(t,x,\alpha)=\|\alpha\|^{2},\quad g(x)=\|x\|^{2},\quad b(t,x,\alpha)=x+\alpha \quad\text{
			for all }(t,x,\alpha)\in [0,T]\times\RR^{d}\times A.$$ 
		In this case, the value function $V$ is characterized by the following HJB equation:
		\begin{equation}
			\begin{aligned}
				-\partial_{t}V(t,x)+\frac{1}{4} \left\|\nabla_{x}V(t,x)\right\|^2-x\cdot  \nabla_{x}V(t,x)&= 0\quad\text{for all }(t,x)\in]0,1[\times\RR^{d}, \\
				V\left(1,x\right)&=\|x\|^2\quad\text{for all }x\in\RR^{d},
			\end{aligned}
			\label{hjb_test} 
		\end{equation}
		whose unique (classical) solution is given by 
		\begin{equation}
			V\left(t, x\right)=\frac{2\|x\|^2}{1+e^{2\left(t-1\right)}}\quad\text{for all }(t,x)\in[0,1]\times\RR^{d}. 
		\end{equation}
		
		{\it Test 1.} We first consider equation~\eqref{hjb_test} in the one-dimensional setting $d=1$. Notice that this simple instance could be solved with more precise and standard discretization methods such as finite differences and semi-Lagrangian schemes. However, the purpose here is to compare ML-based methods having at our disposal an explicit solution. 
		
		We run Algorithm~\ref{alg:IGT} over three rounds using a batch of $M=1000$ points  uniformly distributed in the time-space domain $[0,1] \times [-1,1]$ to implement the initialization step with DGM. We consider a batch of size $S=128$ consisting of initial conditions uniformly sampled in $[-1,1]$ to generate the data set $\mathcal{D}_{\rm OC}$. The computation of the term $\lossHJB$ in~\eqref{eq:loss_total}, during the training process, is performed using a batch of points uniformly distributed in the time-space domain $[0,1] \times [-1,1]$ of the same size $M$ than the one in the initialization step. The neural network $V^{\text{NN}}_{\theta}$, given by~\eqref{para hjb}, is parameterized using $\varphi_1$, defined in~\eqref{eq:parametrizations_phi}, and, in the loss function $\mathrm{Loss}_{\mathrm{val}}(\theta)$ (see~\eqref{eq:loss_total}), we set the penalization parameter to $\lambda_1 = 1$.
		
\begin{figure}
			\centering\includegraphics[width=5cm]{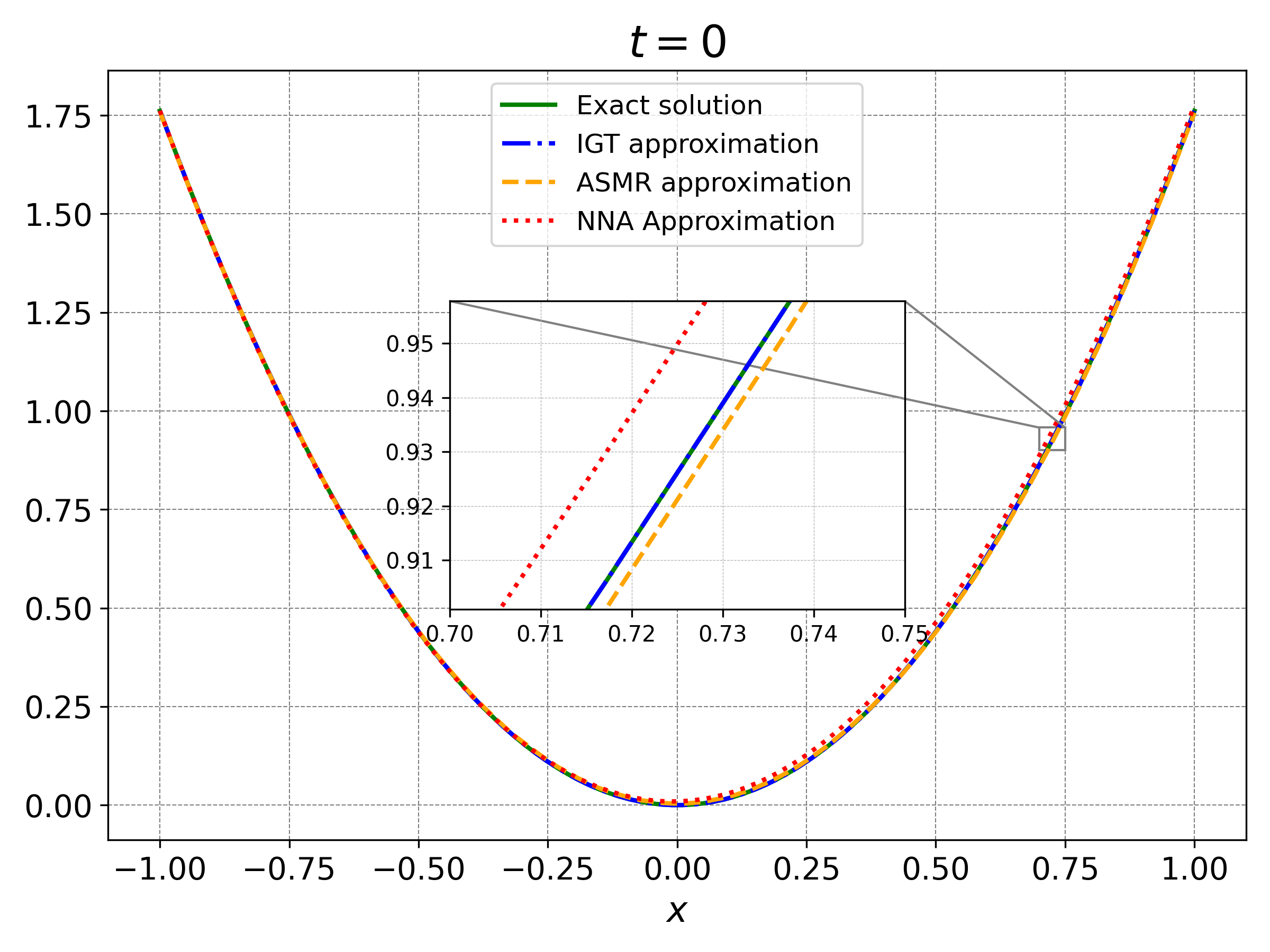}
			\centering\includegraphics[width=5cm]{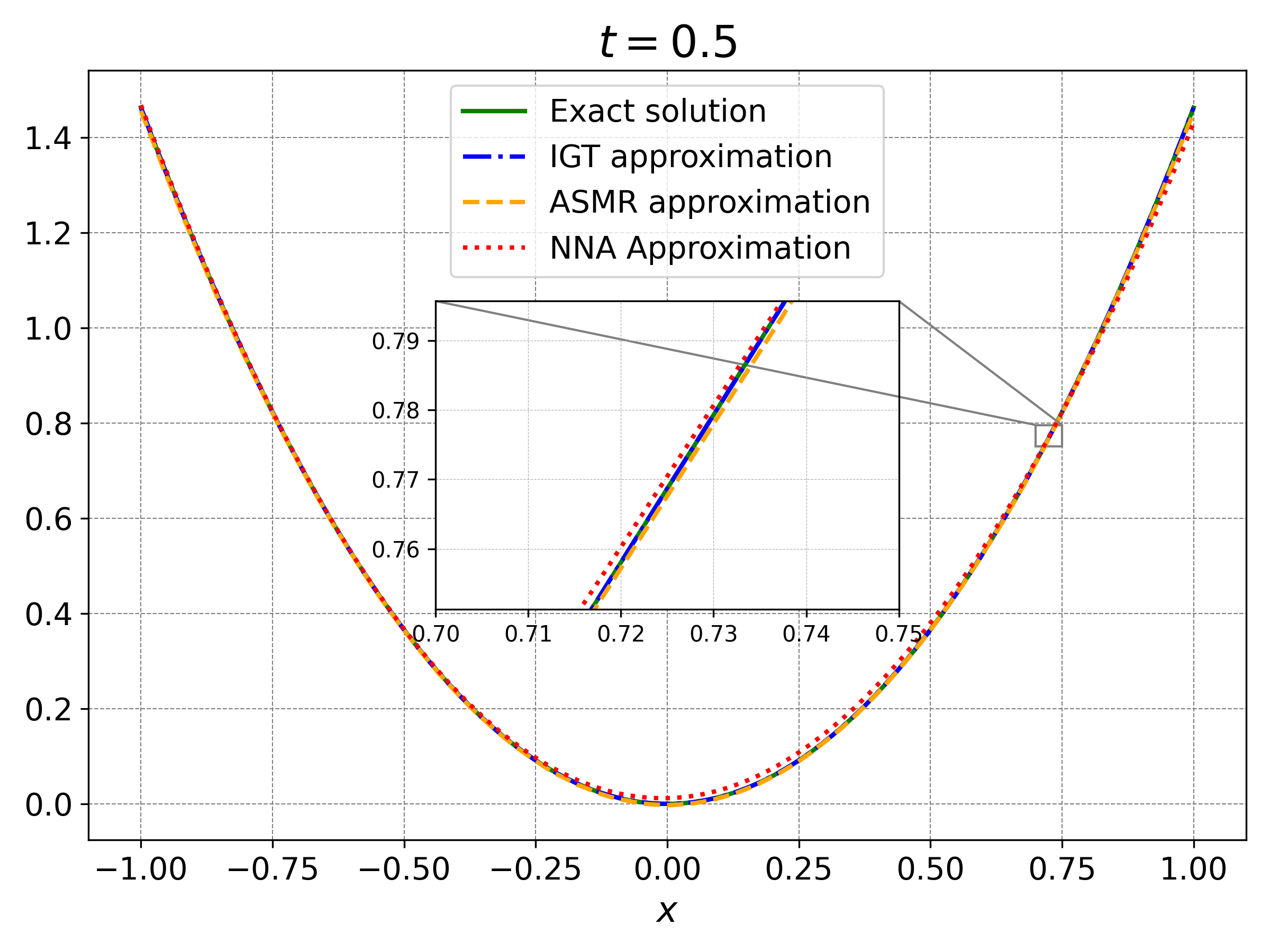}
			\centering\includegraphics[width=5cm]{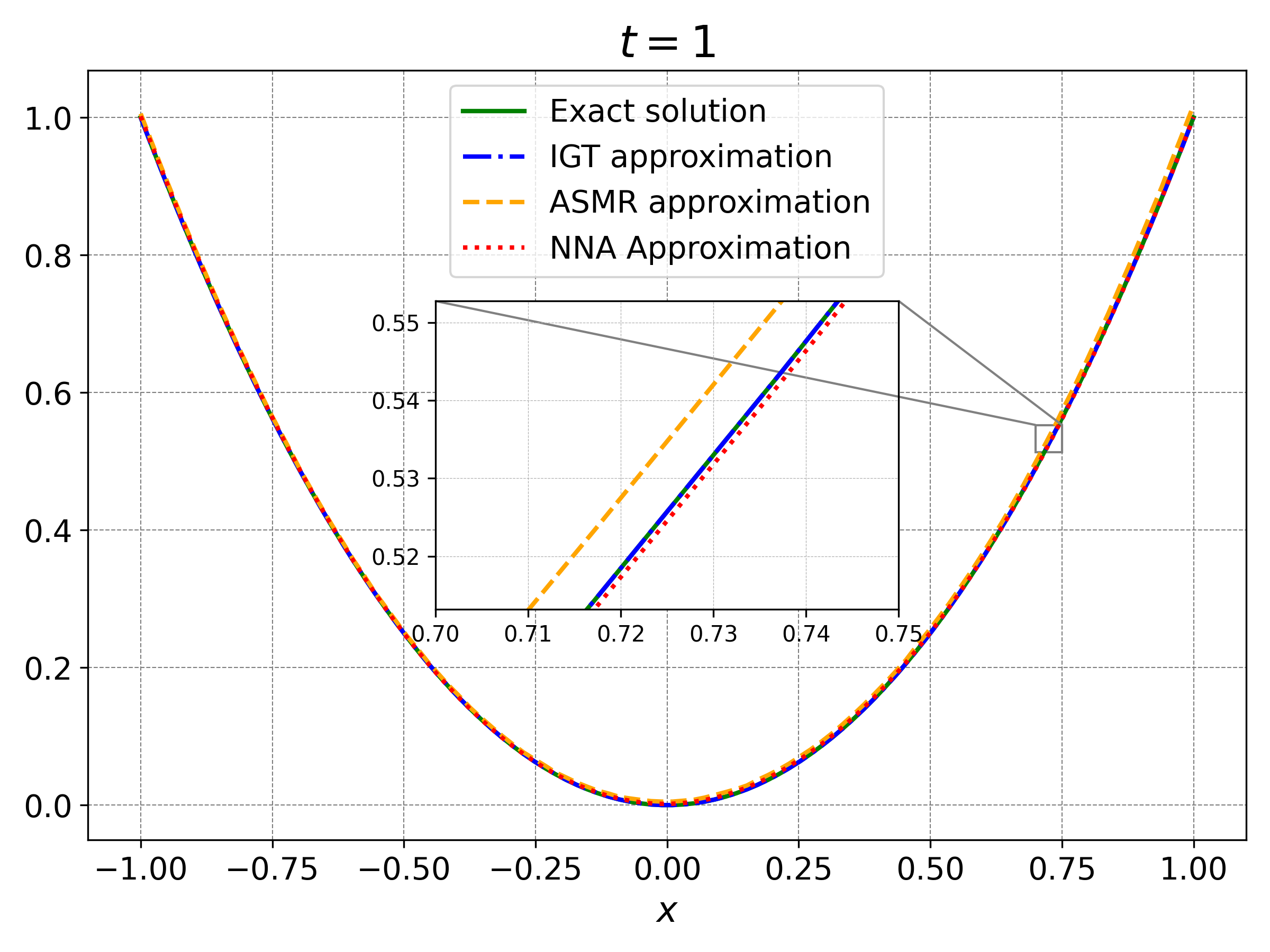}
			\caption{Exact value function and its approximations computed with the IGT, ASMR and NNA methods at three time instances ($d=1$).} 
			\label{fig2}
		\end{figure}

\begin{figure}
			\centering\includegraphics[width=5cm]{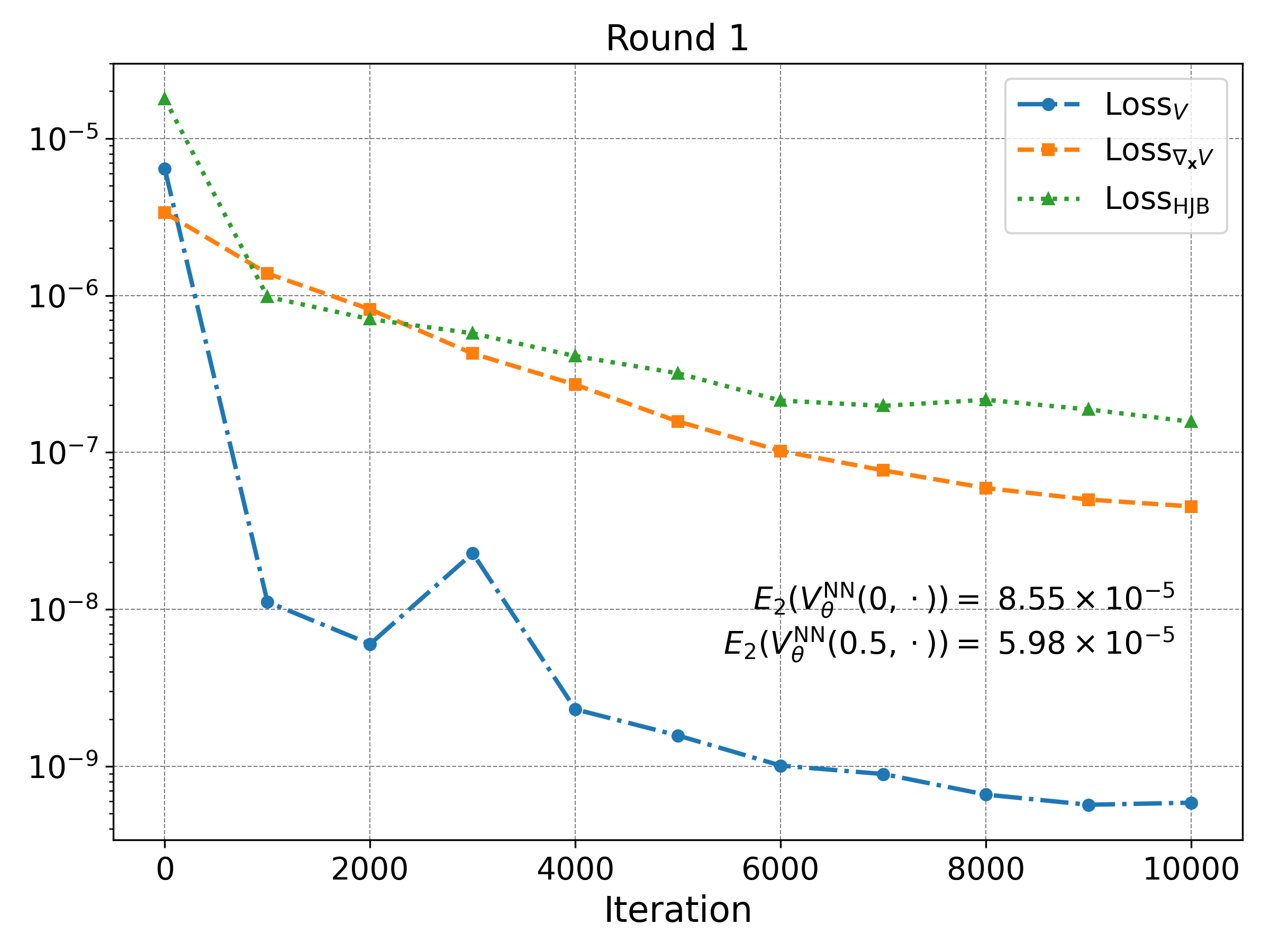}
			\centering\includegraphics[width=5cm]{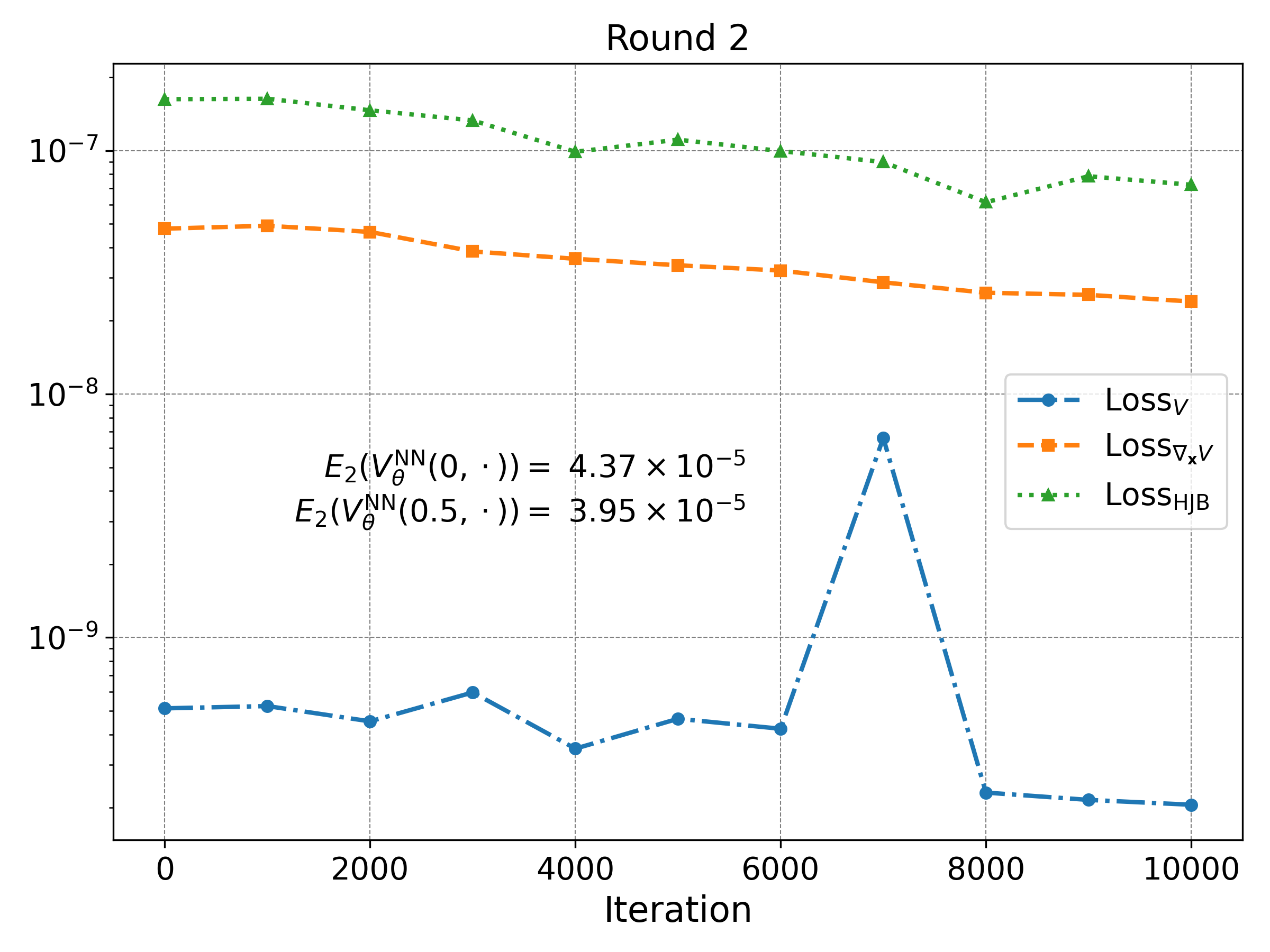}
			\centering\includegraphics[width=5cm]{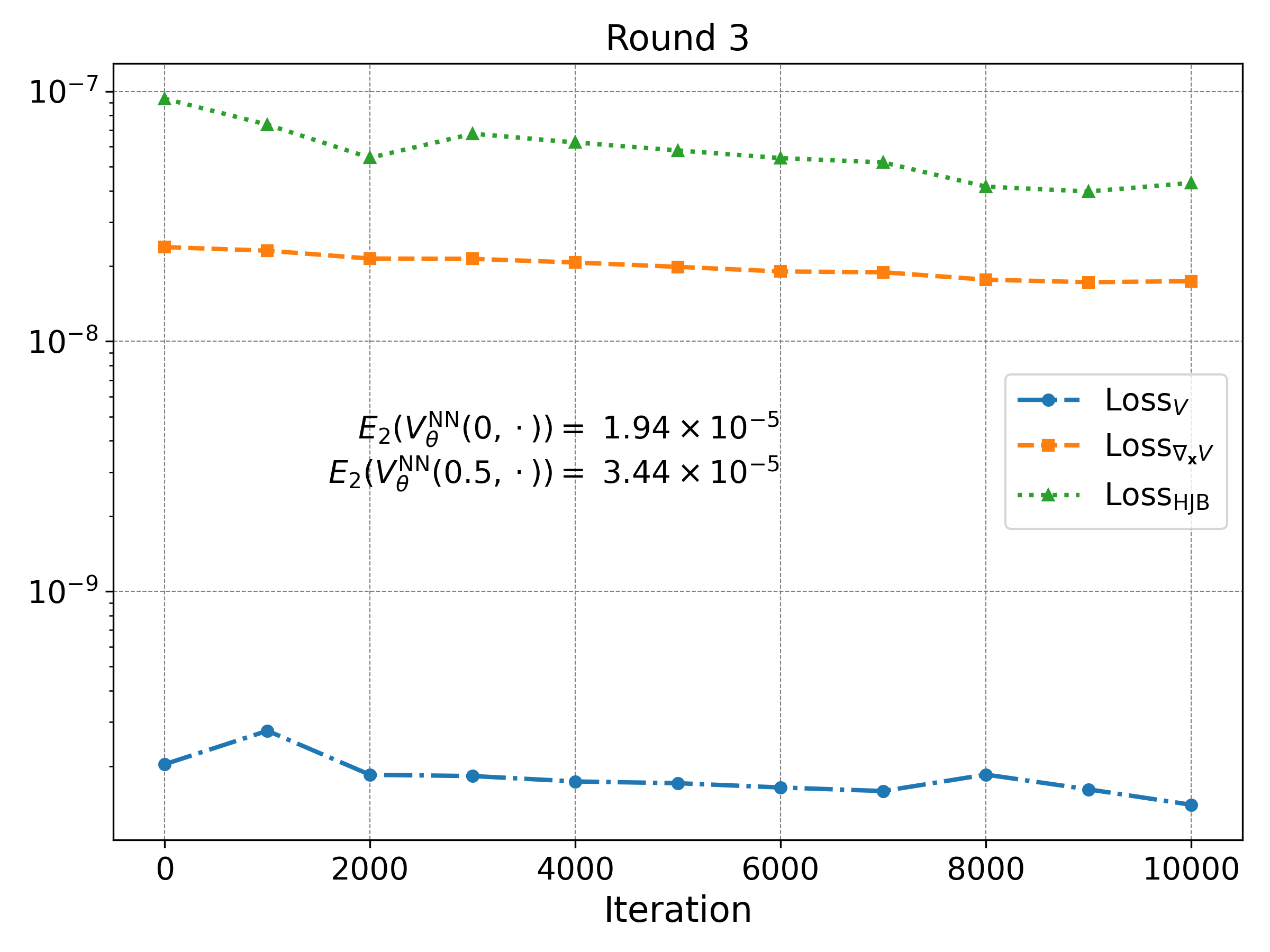}
			\caption{ Residual losses of the training step in the IGT algorithm over three rounds. We also show the relative errors  $E_{2}(V^{\mathrm{NN}}_{\theta}(t,\cdot))$ for $t=0$, $0.5$ at the end of each round.} 
			\label{fig3}
		\end{figure}

		Tables~\ref{tab:comparison_d1} and ~\ref{tab:control_comparison_d1} report the relative errors $E_{2}(V^{\mathrm{NN}}_{\theta}(t,\cdot))$ and  $E_{\infty}(V^{\mathrm{NN}}_{\theta}(t,\cdot))$  for the value function, and $E_{2}(\alpha_{\theta}(t,\cdot))$ and $E_{\infty}(\alpha_{\theta}(t,\cdot))$ for the control, respectively for the IGT, ASMR and NNA methods at different time instances, computed over a uniform spatial grid of 1000 points in the interval $[-1,1]$. Figure~\ref{fig2} shows the corresponding numerical results, illustrating the convergence of the methods at three representative time instants. To monitor the convergence of the IGT method, Figure~\ref{fig3} displays the evolution of the three residual losses defined in~\eqref{eq:loss_V_grad_HJB} over the training iterations at each round. These curves help to assess whether the training process has converged or if further iterations are needed.
		
		\begin{table}[htbp]
			\centering
			\caption{Relative errors of the value function approximation for the IGT, ASMR, and NNA methods at different time instances ($d=1$).}
			\label{tab:comparison_d1}
			\begin{tabular}{|c|ccc|ccc|}
				\hline
				\multirow{2}{*}{\textbf{Time}} 
				& \multicolumn{3}{c|}{$E_{2}$} 
				& \multicolumn{3}{c|}{$E_{\infty}$} \\
				& \textbf{IGT} & \textbf{ASMR} & \textbf{NNA}
				& \textbf{IGT} & \textbf{ASMR} & \textbf{NNA} \\
				\hline
				$t = 0.00$ 
				& $1.94 \times 10^{-5}$ & $7.99 \times 10^{-4}$ & $3.13 \times 10^{-2}$
				& $2.32 \times 10^{-5}$ & $9.12 \times 10^{-4}$ & $2.48 \times 10^{-2}$ \\
				$t = 0.25$ 
				& $1.85 \times 10^{-5}$ & $7.60 \times 10^{-4}$ & $3.54 \times 10^{-2}$
				& $2.14 \times 10^{-5}$ & $2.56 \times 10^{-3}$ & $3.08 \times 10^{-2}$ \\
				$t = 0.50$ 
				& $3.44 \times 10^{-5}$ & $1.17 \times 10^{-3}$ & $4.96 \times 10^{-2}$
				& $9.75 \times 10^{-5}$ & $3.52 \times 10^{-3}$ & $5.23 \times 10^{-2}$ \\
				$t = 0.75$ 
				& $3.35 \times 10^{-5}$ & $3.33 \times 10^{-3}$ & $5.53 \times 10^{-2}$
				& $4.87 \times 10^{-5}$ & $7.71 \times 10^{-3}$ & $5.91 \times 10^{-2}$ \\
				$t = 1.00$ 
				& $0.00$               & $1.07 \times 10^{-2}$ & $1.04 \times 10^{-3}$
				& $0.00$               & $2.15 \times 10^{-2}$ & $1.18 \times 10^{-3}$ \\
				\hline
			\end{tabular}
		\end{table}

		\begin{table}[htbp]
			\centering
			\caption{Relative errors of the feedback control approximation for the IGT, ASMR, and NNA methods at different time instances ($d=1$).}
			\label{tab:control_comparison_d1}
			\begin{tabular}{|c|ccc|ccc|}
				\hline
				\multirow{2}{*}{\textbf{Time}} 
				& \multicolumn{3}{c|}{$E_{2}$} 
				& \multicolumn{3}{c|}{$E_{\infty}$} \\
				& \textbf{IGT} & \textbf{ASMR} & \textbf{NNA}
				& \textbf{IGT} & \textbf{ASMR} & \textbf{NNA} \\
				\hline
				$t = 0.00$ 
				& $1.40 \times 10^{-4}$ & $2.86 \times 10^{-3}$ & $3.71 \times 10^{-2}$
				& $6.42 \times 10^{-4}$ & $1.00 \times 10^{-2}$ & $3.26 \times 10^{-2}$ \\
				$t = 0.25$ 
				& $1.13 \times 10^{-4}$ & $4.86 \times 10^{-3}$ & $4.62 \times 10^{-2}$
				& $4.52 \times 10^{-4}$ & $1.95 \times 10^{-2}$ & $4.23 \times 10^{-2}$ \\
				$t = 0.50$ 
				& $2.10 \times 10^{-4}$ & $5.73 \times 10^{-3}$ & $6.84 \times 10^{-2}$
				& $1.15 \times 10^{-3}$ & $2.16 \times 10^{-2}$ & $7.03 \times 10^{-2}$ \\
				$t = 0.75$ 
				& $2.12 \times 10^{-4}$ & $9.41 \times 10^{-3}$ & $7.30 \times 10^{-2}$
				& $4.65 \times 10^{-4}$ & $2.33 \times 10^{-2}$ & $7.01 \times 10^{-2}$ \\
				$t = 1.00$ 
				& $0.00$               & $3.11 \times 10^{-2}$ & $1.64 \times 10^{-3}$
				& $0.00$               & $6.89 \times 10^{-2}$ & $2.08 \times 10^{-3}$ \\
				\hline
			\end{tabular}
		\end{table}
		
		{\it Test 2.} In this test, we evaluate the performance of the IGT method to solve~\eqref{hjb_test} in several state dimensions $d = 2, 10, 50, 100$, and we compare it with the ASMR and NNA methods.  In all the simulations, we adopt DGM method for the initialization step, and we take $\varphi_1$ as transition function in~\eqref{para hjb} and the batch of points to penalize the residual of the HJB equation and to generate $\mathcal{D}_{\rm OC}$ are sampled uniformly from $[0,1] \times [-1,1]^d$ and $[-1,1]^d$, respectively. For $d = 2$, we consider the same batches sizes $M$ and $S$ than those in the one-dimensional test. We also take $\lambda_1=1$ and we run the algorithm over 3 rounds. In higher dimensions, we run 5 training rounds using batch sizes adapted to the problem complexity. More precisely,  for $d = 10$, we set $M = 2000$, $S = 264$, and $\lambda_1 = 0.5$. For dimensions $d = 50$ and $d = 100$, we increase the batch size to $M = 5000$ while keeping $S = 264$, with penalization parameters $\lambda_1 = 0.1$ and $0.01$, respectively. The corresponding numerical results are detailed in Tables~\ref{tab:comparison_d2_d10_d50_d100_timesteps} and~\ref{tab:comparison_control_d2_d10_d50_d100_timesteps}, where N.C. indicates that the method did not converge within a reasonable number of training iterations.
		
		\begin{table}[htbp]
			\centering
			\caption{Relative errors of the value function approximation for the IGT, ASMR, and NNA methods at different time instances for dimensions $d=2,10,50,100$.}
			\begin{tabular}{|c|c|ccc|ccc|}
				\hline
				\textbf{Dim} & \textbf{Time} 
				& \multicolumn{3}{c|}{$E_2$} 
				& \multicolumn{3}{c|}{$E_\infty$} \\
				& 
				& \textbf{IGT} & \textbf{ASMR} & \textbf{NNA}
				& \textbf{IGT} & \textbf{ASMR} & \textbf{NNA} \\
				\hline
				\multirow{5}{*}{$d = 2$}
				& $t=0.00$ & $1.14\times10^{-4}$ & $1.27\times10^{-3}$ & $6.45\times10^{-2}$
				& $3.83\times10^{-4}$ & $3.15\times10^{-3}$ & $4.89\times10^{-2}$ \\
				& $t=0.25$ & $8.83\times10^{-5}$ & $8.00\times10^{-4}$ & $6.14\times10^{-2}$
				& $2.57\times10^{-4}$ & $4.46\times10^{-3}$ & $4.87\times10^{-2}$ \\
				& $t=0.50$ & $1.02\times10^{-4}$ & $2.00\times10^{-3}$ & $6.88\times10^{-2}$
				& $2.36\times10^{-4}$ & $8.47\times10^{-3}$ & $8.25\times10^{-2}$ \\
				& $t=0.75$ & $1.33\times10^{-4}$ & $2.40\times10^{-3}$ & $7.25\times10^{-2}$
				& $4.67\times10^{-4}$ & $9.59\times10^{-3}$ & $1.01\times10^{-1}$ \\
				& $t=1.00$ & $0.00$              & $8.61\times10^{-3}$ & $7.88\times10^{-4}$
				& $0.00$              & $3.59\times10^{-2}$ & $9.61\times10^{-4}$ \\
				\hline
				
				\multirow{5}{*}{$d = 10$}
				& $t=0.00$ & $1.11\times10^{-3}$ & $1.08\times10^{-3}$ & $1.56\times10^{-1}$
				& $4.59\times10^{-3}$ & $7.88\times10^{-3}$ & $2.36\times10^{-1}$ \\
				& $t=0.25$ & $9.77\times10^{-4}$ & $9.77\times10^{-4}$ & $1.15\times10^{-1}$
				& $4.05\times10^{-3}$ & $6.30\times10^{-3}$ & $1.90\times10^{-1}$ \\
				& $t=0.50$ & $8.78\times10^{-4}$ & $1.19\times10^{-3}$ & $8.67\times10^{-2}$
				& $3.72\times10^{-3}$ & $7.70\times10^{-3}$ & $1.60\times10^{-1}$ \\
				& $t=0.75$ & $6.93\times10^{-4}$ & $1.66\times10^{-3}$ & $5.56\times10^{-2}$
				& $3.44\times10^{-3}$ & $9.50\times10^{-3}$ & $1.08\times10^{-1}$ \\
				& $t=1.00$ & $0.00$              & $2.99\times10^{-3}$ & $5.51\times10^{-3}$
				& $0.00$              & $1.27\times10^{-2}$ & $4.75\times10^{-3}$ \\
				\hline
				
				\multirow{5}{*}{$d = 50$}
				& $t=0.00$ & $6.38\times10^{-3}$ & $1.32\times10^{-2}$ & $3.55\times10^{-1}$
				& $2.78\times10^{-2}$ & $1.03\times10^{-1}$ & $4.63\times10^{-1}$ \\
				& $t=0.25$ & $5.59\times10^{-3}$ & $1.29\times10^{-2}$ & $3.17\times10^{-1}$
				& $2.65\times10^{-2}$ & $1.12\times10^{-1}$ & $4.75\times10^{-1}$ \\
				& $t=0.50$ & $4.67\times10^{-3}$ & $1.24\times10^{-2}$ & $3.16\times10^{-1}$
				& $2.34\times10^{-2}$ & $1.15\times10^{-1}$ & $4.71\times10^{-1}$ \\
				& $t=0.75$ & $3.18\times10^{-3}$ & $1.52\times10^{-2}$ & $2.85\times10^{-1}$
				& $1.62\times10^{-2}$ & $1.18\times10^{-1}$ & $4.14\times10^{-1}$ \\
				& $t=1.00$ & $0.00$              & $4.04\times10^{-2}$ & $1.16\times10^{-1}$
				& $0.00$              & $1.08\times10^{-1}$ & $4.49\times10^{-1}$ \\
				
				\hline
				
				\multirow{5}{*}{$d = 100$}
				& $t=0.00$ & $1.79\times10^{-2}$ & $5.59\times10^{-2}$ & \text{N.C.}
				& $8.02\times10^{-2}$ & $1.04\times10^{-1}$ & \text{N.C.} \\
				& $t=0.25$ & $1.55\times10^{-2}$ & $2.71\times10^{-2}$ & \text{N.C.}
				& $6.94\times10^{-2}$ & $1.34\times10^{-1}$ & \text{N.C.} \\
				& $t=0.50$ & $1.26\times10^{-2}$ & $2.77\times10^{-2}$ & \text{N.C.}
				& $5.60\times10^{-2}$ & $1.49\times10^{-1}$ & \text{N.C.} \\
				& $t=0.75$ & $8.51\times10^{-3}$ & $6.16\times10^{-2}$ & \text{N.C.}
				& $3.68\times10^{-2}$ & $1.61\times10^{-1}$ & \text{N.C.} \\
				& $t=1.00$ & $0.00$                 & $1.47\times10^{-1}$ & \text{N.C.}
				& $0.00$                 & $3.15\times10^{-1}$ & \text{N.C.} \\
				\hline
				
			\end{tabular}
			\label{tab:comparison_d2_d10_d50_d100_timesteps}
		\end{table}
		
		\begin{table}[htbp]
			\centering
			\caption{Relative errors of the feedback control approximation for the IGT, ASMR, and NNA methods at different time instances for dimensions $d = 2$, $10$, $50$, and $100$.}
			\begin{tabular}{|c|c|ccc|ccc|}
				\hline
				\textbf{Dim} & \textbf{Time} 
				& \textbf{IGT} & \textbf{ASMR} & \textbf{NNA}
				& \textbf{IGT} & \textbf{ASMR} & \textbf{NNA} \\
				\hline
				\multicolumn{2}{|c|}{} & \multicolumn{3}{c|}{$E_2$} & \multicolumn{3}{c|}{$E_\infty$} \\
				\hline
				\multirow{5}{*}{$d = 2$}
				& $t=0.00$ & $5.71\times10^{-4}$ & $4.10\times10^{-3}$ & $6.19\times10^{-2}$
				& $6.27\times10^{-3}$ & $2.79\times10^{-2}$ & $6.67\times10^{-2}$ \\
				& $t=0.25$ & $4.67\times10^{-4}$ & $5.03\times10^{-3}$ & $8.44\times10^{-2}$
				& $4.06\times10^{-3}$ & $3.04\times10^{-2}$ & $1.10\times10^{-1}$ \\
				& $t=0.50$ & $5.89\times10^{-4}$ & $7.81\times10^{-3}$ & $1.20\times10^{-1}$
				& $5.15\times10^{-3}$ & $3.81\times10^{-2}$ & $1.61\times10^{-1}$ \\
				& $t=0.75$ & $7.77\times10^{-4}$ & $8.36\times10^{-3}$ & $1.37\times10^{-1}$
				& $4.91\times10^{-3}$ & $4.09\times10^{-2}$ & $1.64\times10^{-1}$ \\
				& $t=1.00$ & $0.00$ & $2.06\times10^{-2}$ & $2.05\times10^{-3}$
				& $0.00$ & $8.64\times10^{-2}$ & $2.85\times10^{-3}$ \\
				\hline
				
				\multirow{5}{*}{$d = 10$}
				& $t=0.00$ & $5.35\times10^{-3}$ & $5.08\times10^{-3}$ & $4.24\times10^{-1}$
				& $3.42\times10^{-2}$ & $3.76\times10^{-2}$ & $4.39\times10^{-1}$ \\
				& $t=0.25$ & $4.88\times10^{-3}$ & $4.75\times10^{-3}$ & $3.84\times10^{-1}$
				& $3.23\times10^{-2}$ & $3.87\times10^{-2}$ & $3.99\times10^{-1}$ \\
				& $t=0.50$ & $4.67\times10^{-3}$ & $5.28\times10^{-3}$ & $3.15\times10^{-1}$
				& $2.71\times10^{-2}$ & $3.93\times10^{-2}$ & $3.27\times10^{-1}$ \\
				& $t=0.75$ & $3.69\times10^{-3}$ & $6.67\times10^{-3}$ & $1.98\times10^{-1}$
				& $2.13\times10^{-2}$ & $4.34\times10^{-2}$ & $2.04\times10^{-1}$ \\
				& $t=1.00$ & $0.00$ & $1.01\times10^{-2}$ & $1.82\times10^{-3}$
				& $0.00$ & $5.66\times10^{-2}$ & $1.34\times10^{-2}$ \\
				\hline
				
				\multirow{5}{*}{$d = 50$}
				& $t=0.00$ & $4.00\times10^{-2}$ & $8.68\times10^{-2}$ & $5.38\times10^{-1}$
				& $1.18\times10^{-1}$ & $3.26\times10^{-1}$ & $5.40\times10^{-1}$ \\
				& $t=0.25$ & $3.54\times10^{-2}$ & $8.17\times10^{-2}$ & $5.96\times10^{-1}$
				& $1.06\times10^{-1}$ & $3.27\times10^{-1}$ & $5.81\times10^{-1}$ \\
				& $t=0.50$ & $2.98\times10^{-2}$ & $8.42\times10^{-2}$ & $6.43\times10^{-1}$
				& $8.78\times10^{-2}$ & $3.33\times10^{-1}$ & $6.08\times10^{-1}$ \\
				& $t=0.75$ & $2.02\times10^{-2}$ & $9.45\times10^{-2}$ & $4.49\times10^{-1}$
				& $5.71\times10^{-2}$ & $3.35\times10^{-1}$ & $5.08\times10^{-1}$ \\
				& $t=1.00$ & $0.00$ & $1.43\times10^{-1}$ & $3.71\times10^{-1}$
				& $0.00$ & $3.30\times10^{-1}$ & $7.24\times10^{-1}$ \\
				\hline
				
				\multirow{5}{*}{$d = 100$}
				& $t=0.00$ & $1.47\times10^{-1}$ & $1.05\times10^{-1}$ & \text{N.C.}
				& $5.77\times10^{-1}$ & $3.80\times10^{-1}$ & \text{N.C.} \\
				& $t=0.25$ & $1.29\times10^{-1}$ & $1.18\times10^{-1}$ & \text{N.C.}
				& $4.50\times10^{-1}$ & $3.84\times10^{-1}$ & \text{N.C.} \\
				& $t=0.50$ & $1.08\times10^{-1}$ & $1.91\times10^{-1}$ & \text{N.C.}
				& $3.21\times10^{-1}$ & $4.02\times10^{-1}$ & \text{N.C.} \\
				& $t=0.75$ & $7.35\times10^{-2}$ & $3.40\times10^{-1}$ & \text{N.C.}
				& $2.43\times10^{-1}$ & $4.68\times10^{-1}$ & \text{N.C.} \\
				& $t=1.00$ & $0.00$ & $5.94\times10^{-1}$ & \text{N.C.}
				& $0.00$ & $7.23\times10^{-1}$ & \text{N.C.} \\
				\hline
			\end{tabular}
			\label{tab:comparison_control_d2_d10_d50_d100_timesteps}
		\end{table}
		
		\bigskip
		{\it Non-differentiable value function.} To further assess the robustness of the IGT methodology, we consider a second benchmark problem whose value function is not differentiable.  Despite the lack of smoothness, our approach maintains strong performance in this framework. This setting stands in sharp contrast with the smooth linear-quadratic example and provides a far more stringent test for learning-based solvers.
		
		We consider the family of optimal control problems~\eqref{def:oc_problem} with $T=1$, $A=\{\alpha\in\RR^{d}\,:\,\|\alpha\|\leq 1\}$, and 
		$$
		\ell(t,x,\alpha)=\|\alpha\|^{2},\quad g(x)=-\|x\|^{2},\quad b(t,x,\alpha)=\alpha \quad\text{
			for all }(t,x,\alpha)\in [0,T]\times\RR^{d}\times A.$$ 
		
		The associated value function is given by 
		$$
		V(t,x)=\inf\left\{\int_t^1 \|\alpha(s)\|^2\,\dd s - \|\x(1)\|^2 \,\Big|\, \alpha\in L^\infty([0,T];\RR^{d}),\,\|\alpha\|_{\infty}\leq 1 \right\},
		$$
		where $\x$ solves the controlled dynamics with initial condition $\x(t)=x$. Notice that, by Jensen's inequality,
		\begin{equation}
		\label{eq:ex_nonlineaire}
		\begin{aligned}
		V(t,x)&\geq \inf_{\alpha}\left\{\frac{1}{1-t}\left\|\int_{t}^{1}\alpha(s) \dd s\right\|^2-\left\|x+\int_{t}^{1}\alpha(s) \dd s\right\|^{2}\,\Big|\,\, \alpha\in L^\infty([0,T];\RR^{d}),\,\|\alpha\|_{\infty}\leq 1 \right\}\\
		&=\inf\left\{\frac{\|y\|^2}{1-t}-\|x+y\|^{2}\,\Big|\,\alpha\in L^\infty([0,T];\RR^{d}),\,\|\alpha\|_{\infty}\leq 1,\,y=\int_{t}^{1}\alpha(s)\dd s\right\}\\
		&=\inf\left\{\frac{\|y\|^2}{1-t}-\|x+y\|^{2}\,\Big|\,y\in\RR^{d},\,\|y\|\leq 1-t\right\}
        \end{aligned}
		\end{equation}
	Given a minimizer $y^*$ in the last infimum, the constant control $\alpha^*(s)\equiv y^*/(1-t)$ satisfies 
		$$
		\int_t^1 \|\alpha^*(s)\|^2\,\dd s =\frac{1}{1-t}\left\|\int_{t}^{1}\alpha^*(s)\dd s\right\|^2
		$$
		and hence $\alpha^*$ is optimal for $V(t,x)$ and equality holds in~\eqref{eq:ex_nonlineaire}. A direct computation using the last infimum in~\eqref{eq:ex_nonlineaire} then yields
		$$
		V(t,x)=\begin{cases}
			-\dfrac{\|x\|^2}{t}, & t\in ]0,1],\ \|x\|\le t,\\[8pt]
			-\|x\|^2-2(1-t)\|x\|+(1-t)t, & t\in ]0,1],\ \|x\|>t,\\[8pt]
			-\|x\|^2-2\|x\|, & t=0.
		\end{cases}
		$$
		For every $t\in ]0,1]$, the function $V(t,\cdot)$ belongs to $C^\infty(\mathbb{R}^d)$, whereas $V(0,\cdot)$ fails to be differentiable at $x=0$.
		
		To compare numerical performance, we evaluate not only the errors in approximating the value function, but also those incurred in approximating optimal feedback controls. 
		Using the explicit expression of $V$, we obtain 
		$$
		\alpha^*(t,x)=
		\begin{cases}
			-\dfrac{x}{t}, & t\in ]0,1],\ \|x\|\le t,\\[8pt]
			-\dfrac{x}{\|x\|}, & t\in ]0,1],\ \|x\|> t,\\[8pt]
			-\dfrac{x}{\|x\|}, & t=0,\ \|x\|\neq 0,
		\end{cases}
		$$
		while at $(t,x)=(0,0)$ the feedback is not uniquely defined. The lack of uniqueness of the feedback at $(0,0)$ is consistent with the non-differentiability of $V(0,\cdot)$ at the origin.

\begin{figure}
			\centering\includegraphics[width=10cm]{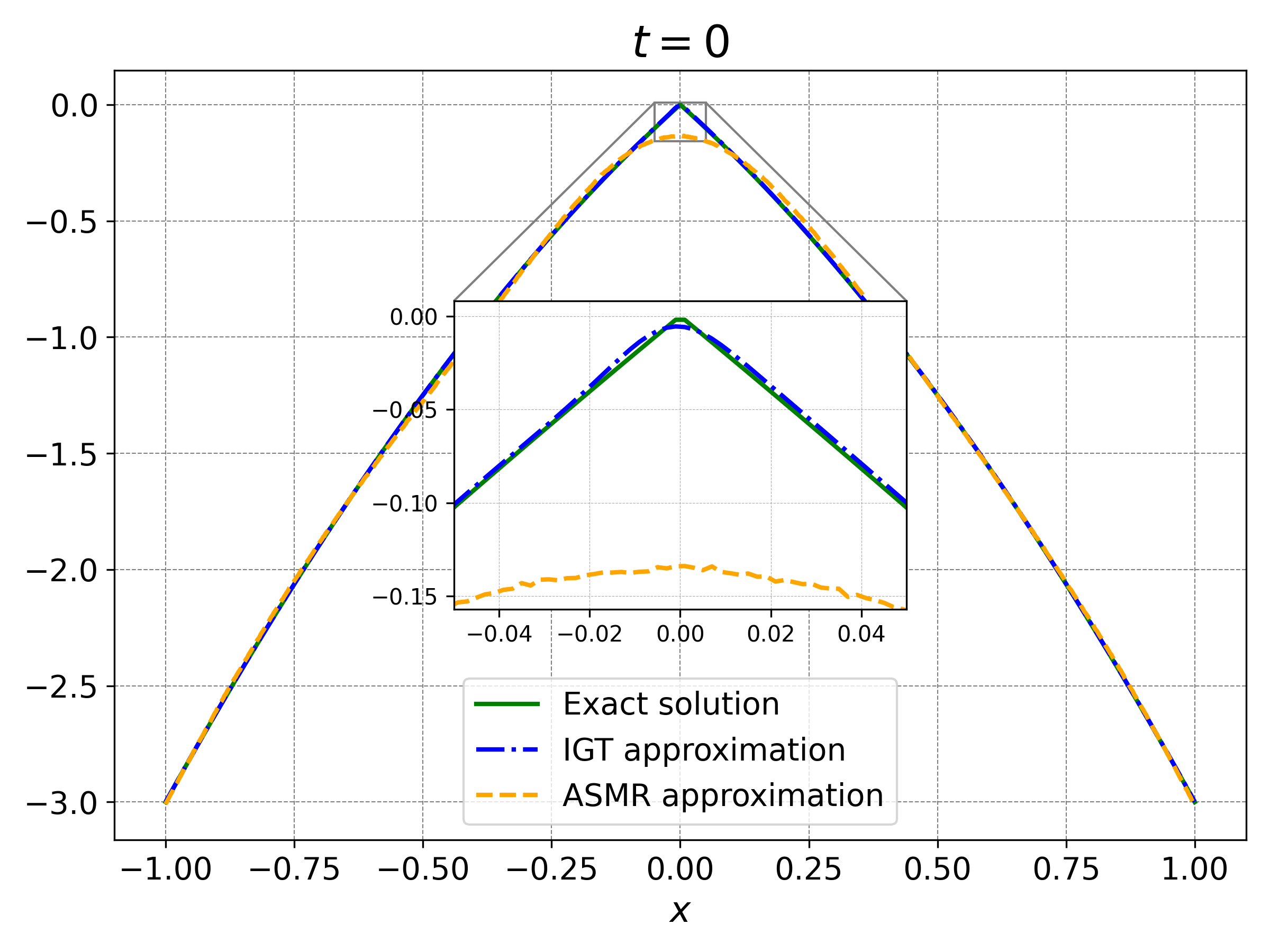}
			
			\caption{Exact value function and its approximation computed with the IGT and ASMR methods at $t=0$.} 
			\label{fig_NN}
		\end{figure}

		In the numerical experiments, we 
		evaluate the IGT and ASMR methods across dimensions $d = 1,2,10,50,100$.  
		To ensure a fair comparison with the smooth LQ tests, we deliberately keep the {\it same training protocol and hyperparameters} as in Test~1 and Test~2, 
		modifying only the parametrization of $V^{\mathrm{NN}}_{\theta}$ by replacing $\varphi_{1}$ with $\varphi_{2}$ in~\eqref{para hjb}.  
		The value-function errors are reported in Table~\ref{tab:comparison_d_V}, and the control errors in Table~\ref{tab:comparison_control_u}.  
		Figure~\ref{fig_NN} displays the numerical approximations at time $t=0$, together with a zoom at the nondifferentiable point $x=0$, where the kink in the exact solution is clearly visible. Notice that we have excluded the NNA method from this experiment because it fails to converge in this nonsmooth problem.
		
		\begin{table}[htbp]
			\centering
			\caption{Relative errors of the value function approximation for the IGT and ASMR methods across dimensions $d = 1,2,10,50,100$.}
			\label{tab:comparison_d_V}
			\begin{tabular}{|c|c|cc|cc|}
				\hline
				\textbf{Dim} & \textbf{Time} 
				& \multicolumn{2}{c|}{$E_2$} 
				& \multicolumn{2}{c|}{$E_\infty$} \\
				& 
				& \textbf{IGT} & \textbf{ASMR}
				& \textbf{IGT} & \textbf{ASMR} \\
				\hline
				
				\multirow{5}{*}{$d=1$}
				& $t=0.00$ & $8.52\times10^{-4}$ & $1.73\times10^{-2}$ & $1.46\times10^{-3}$ & $4.40\times10^{-2}$ \\
				& $t=0.25$ & $8.68\times10^{-4}$ & $3.86\times10^{-2}$ & $8.32\times10^{-4}$ & $4.92\times10^{-2}$ \\
				& $t=0.50$ & $8.36\times10^{-4}$ & $1.46\times10^{-1}$ & $7.25\times10^{-4}$ & $1.59\times10^{-1}$ \\
				& $t=0.75$ & $6.29\times10^{-4}$ & $3.19\times10^{-1}$ & $6.18\times10^{-4}$ & $2.97\times10^{-1}$ \\
				& $t=1.00$ & $0.00$                 & $5.62\times10^{-1}$ & $0.00$                 & $4.60\times10^{-1}$ \\
				\hline
				
				\multirow{5}{*}{$d=2$}
				& $t=0.00$ & $8.41\times10^{-4}$ & $5.67\times10^{-2}$ & $4.68\times10^{-3}$ & $1.54\times10^{-1}$ \\
				& $t=0.25$ & $8.16\times10^{-4}$ & $5.64\times10^{-2}$ & $3.10\times10^{-3}$ & $1.02\times10^{-1}$ \\
				& $t=0.50$ & $8.93\times10^{-4}$ & $5.03\times10^{-2}$ & $1.98\times10^{-3}$ & $6.33\times10^{-2}$ \\
				& $t=0.75$ & $8.53\times10^{-4}$ & $3.54\times10^{-2}$ & $1.29\times10^{-3}$ & $3.37\times10^{-2}$ \\
				& $t=1.00$ & $0.00$                 & $4.99\times10^{-2}$ & $0.00$                 & $4.01\times10^{-2}$ \\
				\hline
				
				\multirow{5}{*}{$d=10$}
				& $t=0.00$ & $4.67\times10^{-3}$ & $8.95\times10^{-3}$ & $2.96\times10^{-2}$ & $6.67\times10^{-2}$ \\
				& $t=0.25$ & $3.72\times10^{-3}$ & $8.39\times10^{-3}$ & $1.84\times10^{-2}$ & $4.95\times10^{-2}$ \\
				& $t=0.50$ & $2.96\times10^{-3}$ & $1.13\times10^{-2}$ & $1.57\times10^{-2}$ & $3.13\times10^{-2}$ \\
				& $t=0.75$ & $1.95\times10^{-3}$ & $1.90\times10^{-2}$ & $9.76\times10^{-3}$ & $2.63\times10^{-2}$ \\
				& $t=1.00$ & $0.00$                 & $4.38\times10^{-2}$ & $0.00$                 & $4.50\times10^{-2}$ \\
				\hline
				
				\multirow{5}{*}{$d=50$}
				& $t=0.00$ & $1.13\times10^{-2}$ & $1.10\times10^{-2}$ & $7.33\times10^{-2}$ & $5.96\times10^{-2}$ \\
				& $t=0.25$ & $8.94\times10^{-3}$ & $1.17\times10^{-2}$ & $5.75\times10^{-2}$ & $6.99\times10^{-2}$ \\
				& $t=0.50$ & $6.44\times10^{-3}$ & $2.01\times10^{-2}$ & $4.14\times10^{-2}$ & $7.89\times10^{-2}$ \\
				& $t=0.75$ & $3.54\times10^{-3}$ & $2.42\times10^{-2}$ & $2.25\times10^{-2}$ & $8.27\times10^{-2}$ \\
				& $t=1.00$ & $0.00$                 & $2.28\times10^{-2}$ & $0.00$                 & $8.28\times10^{-2}$ \\
				\hline
				
				\multirow{5}{*}{$d=100$}
				& $t=0.00$ & $2.04\times10^{-2}$ & $1.24\times10^{-2}$ & $9.50\times10^{-2}$ & $1.12\times10^{-1}$ \\
				& $t=0.25$ & $1.60\times10^{-2}$ & $1.47\times10^{-2}$ & $7.15\times10^{-2}$ & $1.18\times10^{-1}$ \\
				& $t=0.50$ & $1.13\times10^{-2}$ & $1.83\times10^{-2}$ & $4.89\times10^{-2}$ & $1.27\times10^{-1}$ \\
				& $t=0.75$ & $5.98\times10^{-3}$ & $2.67\times10^{-2}$ & $2.56\times10^{-2}$ & $1.41\times10^{-1}$ \\
				& $t=1.00$ & $0.00$                 & $4.09\times10^{-2}$ & $0.00$                 & $1.58\times10^{-1}$ \\
				\hline
				
			\end{tabular}
		\end{table}

		\begin{table}[htbp]
			\centering
			\caption{Relative errors of the feedback control approximation for the IGT and ASMR methods across dimensions $d=1,2,10,50,100$.}
			\label{tab:comparison_control_u}
			\begin{tabular}{|c|c|cc|cc|}
				\hline
				\textbf{Dim} & \textbf{Time} 
				& \multicolumn{2}{c|}{$E_2$}
				& \multicolumn{2}{c|}{$E_\infty$} \\
				& & \textbf{IGT} & \textbf{ASMR} & \textbf{IGT} & \textbf{ASMR} \\
				\hline
				
				\multirow{5}{*}{$d=1$}
				& $t=0.00$ & $5.34\times10^{-2}$ & $2.12\times10^{-1}$ & $9.55\times10^{-1}$ & $1.02\times10^{0}$ \\
				& $t=0.25$ & $4.78\times10^{-3}$ & $8.68\times10^{-2}$ & $3.04\times10^{-2}$ & $1.29\times10^{1}$ \\
				& $t=0.50$ & $2.22\times10^{-3}$ & $2.55\times10^{-1}$ & $1.06\times10^{-2}$ & $2.28\times10^{1}$ \\
				& $t=0.75$ & $3.08\times10^{-3}$ & $3.68\times10^{-1}$ & $5.42\times10^{-3}$ & $1.95\times10^{1}$ \\
				& $t=1.00$ & $0.00$                 & $4.49\times10^{-1}$ & $0.00$                 & $4.56\times10^{0}$ \\
				\hline
				
				\multirow{5}{*}{$d=2$}
				& $t=0.00$ & $1.85\times10^{-2}$ & $1.53\times10^{-1}$ & $8.09\times10^{-1}$ & $9.92\times10^{-1}$ \\
				& $t=0.25$ & $1.07\times10^{-2}$ & $1.41\times10^{-1}$ & $9.14\times10^{-2}$ & $1.07\times10^{0}$ \\
				& $t=0.50$ & $7.11\times10^{-3}$ & $9.57\times10^{-2}$ & $3.84\times10^{-2}$ & $1.47\times10^{0}$ \\
				& $t=0.75$ & $4.89\times10^{-3}$ & $5.01\times10^{-2}$ & $1.54\times10^{-2}$ & $2.31\times10^{0}$ \\
				& $t=1.00$ & $0.00$                 & $3.79\times10^{-2}$ & $0.00$                 & $3.48\times10^{0}$ \\
				\hline
				
				\multirow{5}{*}{$d=10$}
				& $t=0.00$ & $2.28\times10^{-2}$ & $1.73\times10^{-2}$ & $1.53\times10^{-1}$ & $1.00\times10^{-1}$ \\
				& $t=0.25$ & $1.61\times10^{-2}$ & $1.71\times10^{-2}$ & $7.87\times10^{-2}$ & $1.97\times10^{-1}$ \\
				& $t=0.50$ & $1.17\times10^{-2}$ & $2.03\times10^{-2}$ & $1.45\times10^{-1}$ & $2.88\times10^{-1}$ \\
				& $t=0.75$ & $7.11\times10^{-3}$ & $2.55\times10^{-2}$ & $9.67\times10^{-2}$ & $1.62\times10^{-1}$ \\
				& $t=1.00$ & $0.00$                 & $3.32\times10^{-2}$ & $0.00$                 & $1.21\times10^{-1}$ \\
				\hline
				
				\multirow{5}{*}{$d=50$}
				& $t=0.00$ & $8.45\times10^{-2}$ & $6.30\times10^{-2}$ & $3.48\times10^{-1}$ & $2.33\times10^{-1}$ \\
				& $t=0.25$ & $6.53\times10^{-2}$ & $6.23\times10^{-2}$ & $2.64\times10^{-1}$ & $2.23\times10^{-1}$ \\
				& $t=0.50$ & $4.54\times10^{-2}$ & $6.31\times10^{-2}$ & $1.77\times10^{-1}$ & $2.19\times10^{-1}$ \\
				& $t=0.75$ & $2.39\times10^{-2}$ & $6.51\times10^{-2}$ & $8.82\times10^{-2}$ & $2.18\times10^{-1}$ \\
				& $t=1.00$ & $0.00$                 & $6.84\times10^{-2}$ & $0.00$                 & $2.19\times10^{-1}$ \\
				\hline
				
				\multirow{5}{*}{$d=100$}
				& $t=0.00$ & $1.99\times10^{-1}$ & $1.24\times10^{-1}$ & $5.47\times10^{-1}$ & $4.42\times10^{-1}$ \\
				& $t=0.25$ & $1.53\times10^{-1}$ & $1.24\times10^{-1}$ & $4.15\times10^{-1}$ & $4.42\times10^{-1}$ \\
				& $t=0.50$ & $1.05\times10^{-1}$ & $1.24\times10^{-1}$ & $2.77\times10^{-1}$ & $4.42\times10^{-1}$ \\
				& $t=0.75$ & $5.44\times10^{-2}$ & $1.25\times10^{-1}$ & $1.36\times10^{-1}$ & $4.42\times10^{-1}$ \\
				& $t=1.00$ & $0.00$                 & $1.25\times10^{-1}$ & $0.00$                 & $4.42\times10^{-1}$ \\
				\hline
				
			\end{tabular}
		\end{table}

		\subsubsection{{\bf Quadcopter}}
		We illustrate the practical utility of the IGT method by applying it to a well-known real-world benchmark problem, treated in \cite{lin2018splitting, lin2021alternating, onken2022neural}. Specifically, we consider a quadcopter, an aerial vehicle with four rotary wings similar to consumer drones, with dynamics modelled by
		\begin{equation}
			\left\{\begin{array}{l}\ddot{x}=\frac{u}{m}(\sin (\xi) \sin (\phi)+\cos (\xi) \sin (\theta) \cos (\phi)) \\ \ddot{y}=\frac{u}{m}(-\cos (\xi) \sin (\phi)+\sin (\xi) \sin (\theta) \cos (\phi)) \\ \ddot{z}=\frac{u}{m} \cos (\theta) \cos (\phi)-g \\ \ddot{\xi}=\tau_\xi \\ \ddot{\theta}=\tau_\theta \\ 
				\ddot{\phi}=\tau_{\phi},
			\end{array}\right.   
			\label{eq:quadcopter_dynamics}
		\end{equation}
		where $(x,y,z)$ denotes the spatial position of the quadcopter and $\xi$, $\theta$, and $\phi$ are the angular orientations of yaw, pitch, and roll, respectively. The constant $m$ is the mass, which is set to $1$ $(\mathrm{kg})$, and $g=9.81$ $(\mathrm{m}/\mathrm{s}^{2})$ is the standard gravity constant. The system is controlled by the main thrust $u$ and the yaw, pitch, and roll torques denoted by $\tau_\xi$, $\tau_\theta$, and $\tau_\phi$, respectively.
		
\begin{figure}[!htbp]
			\centering\includegraphics[width=10cm]{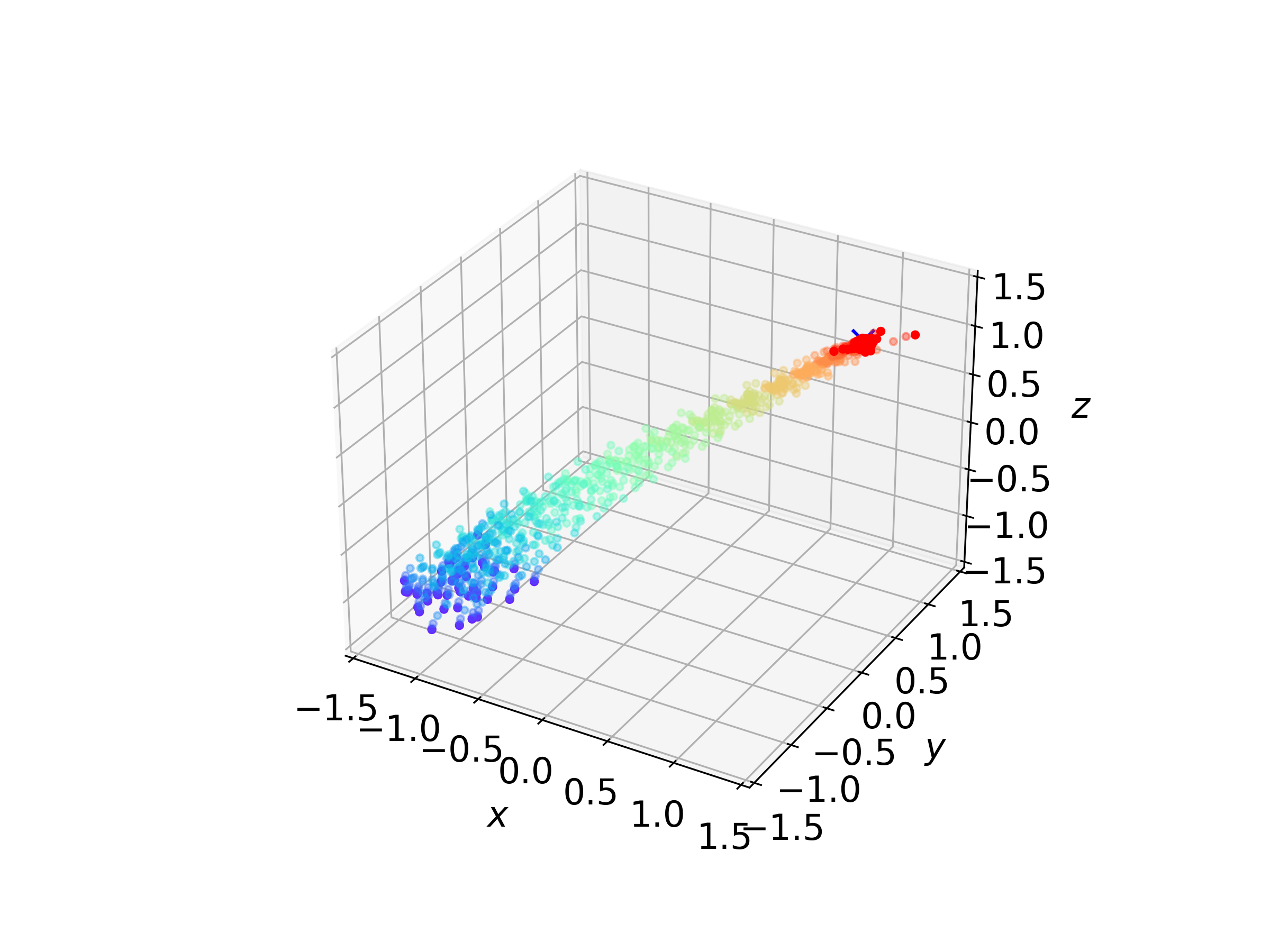}
			\centering\includegraphics[width=7cm]{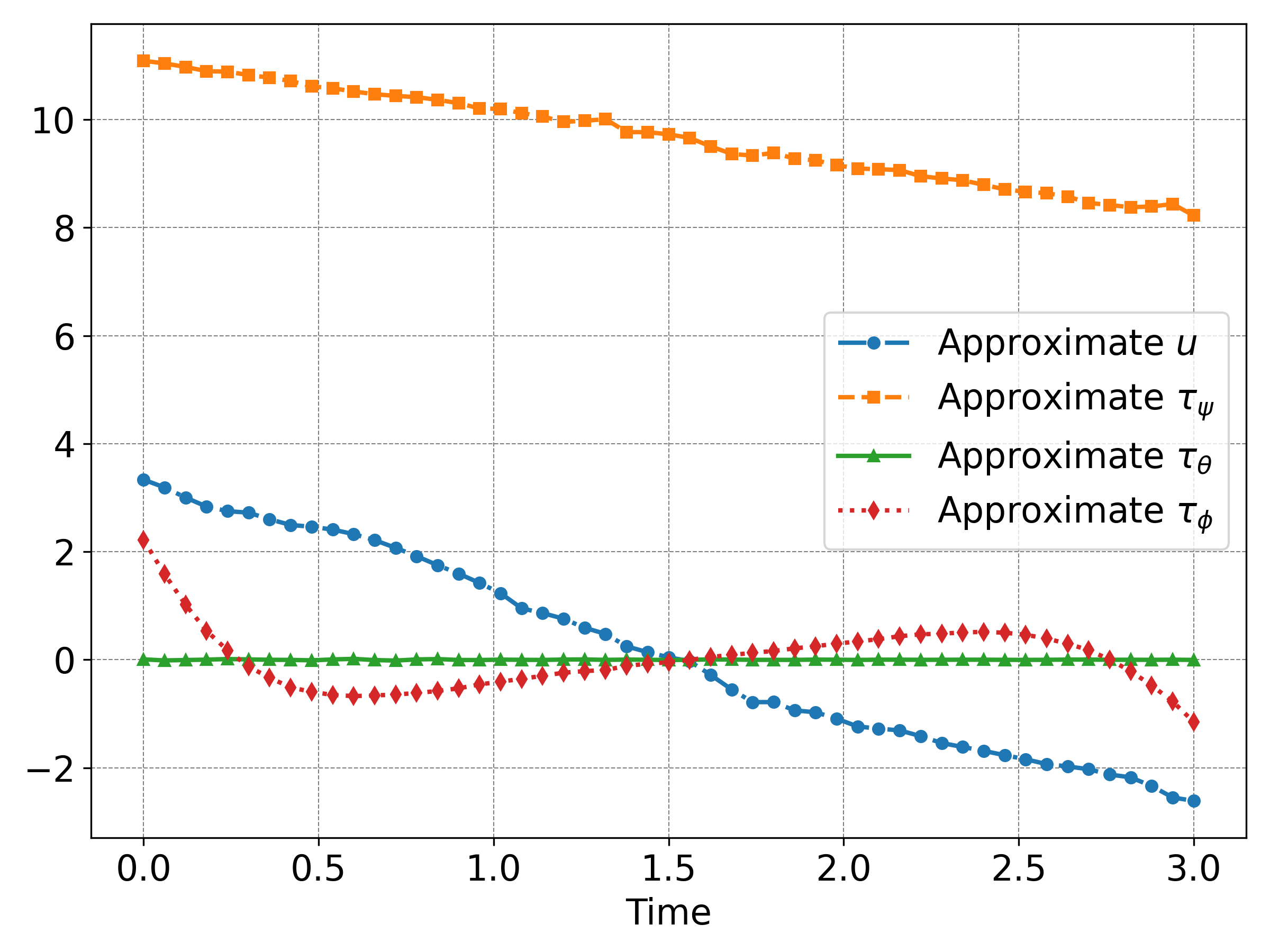}
			\centering\includegraphics[width=7cm]{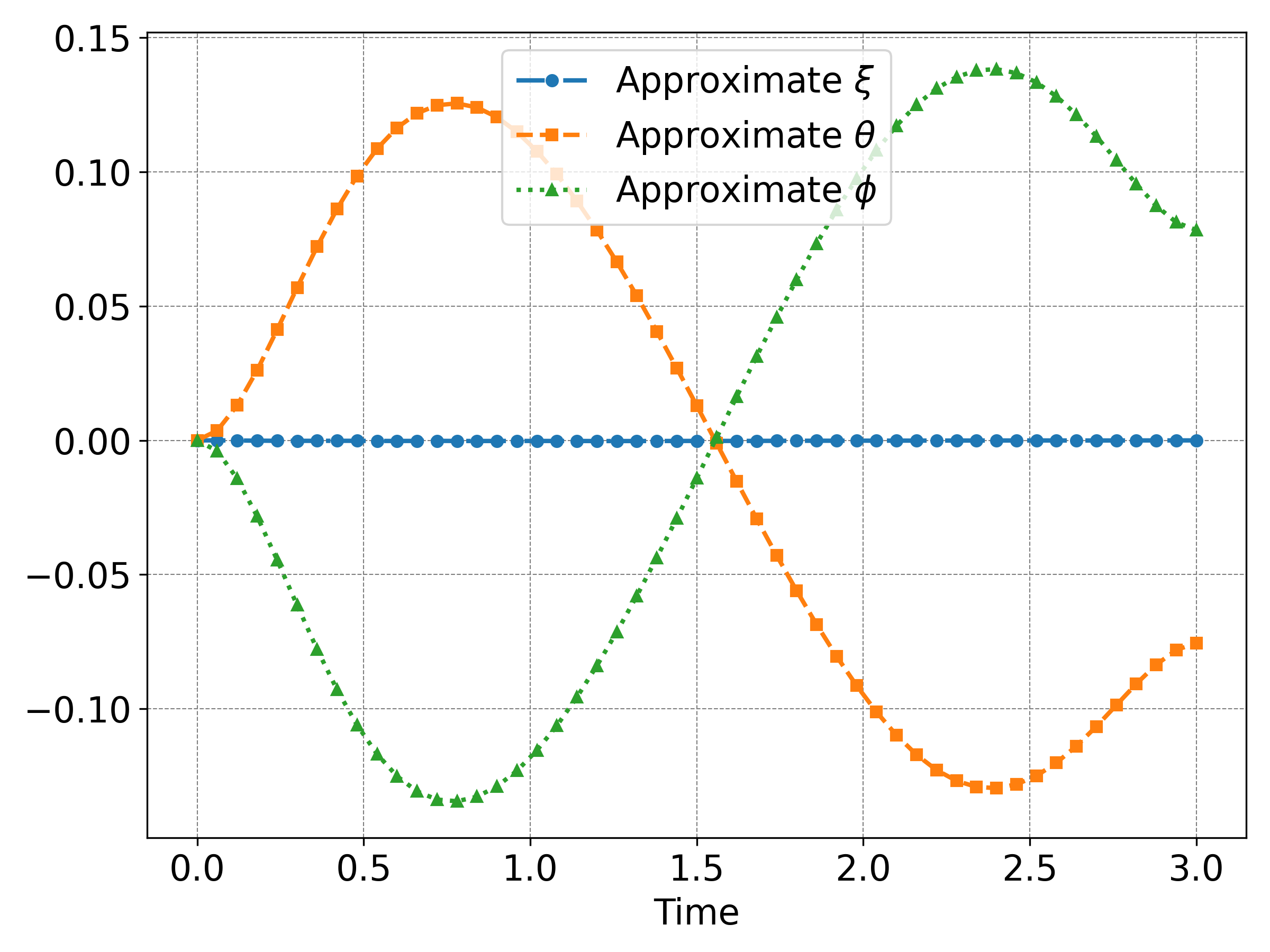}
			\centering\includegraphics[width=7cm]{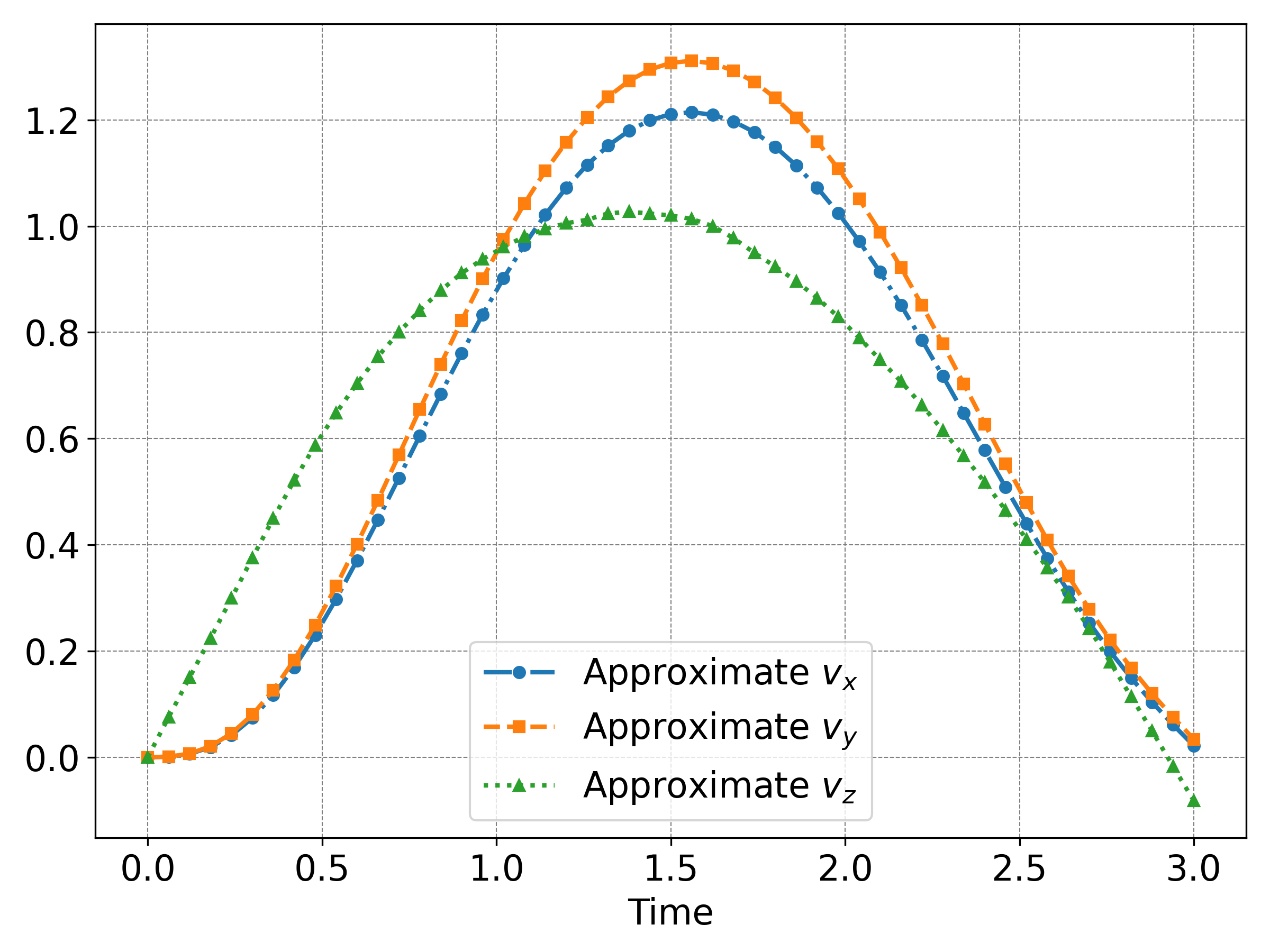}
			\centering\includegraphics[width=7cm]{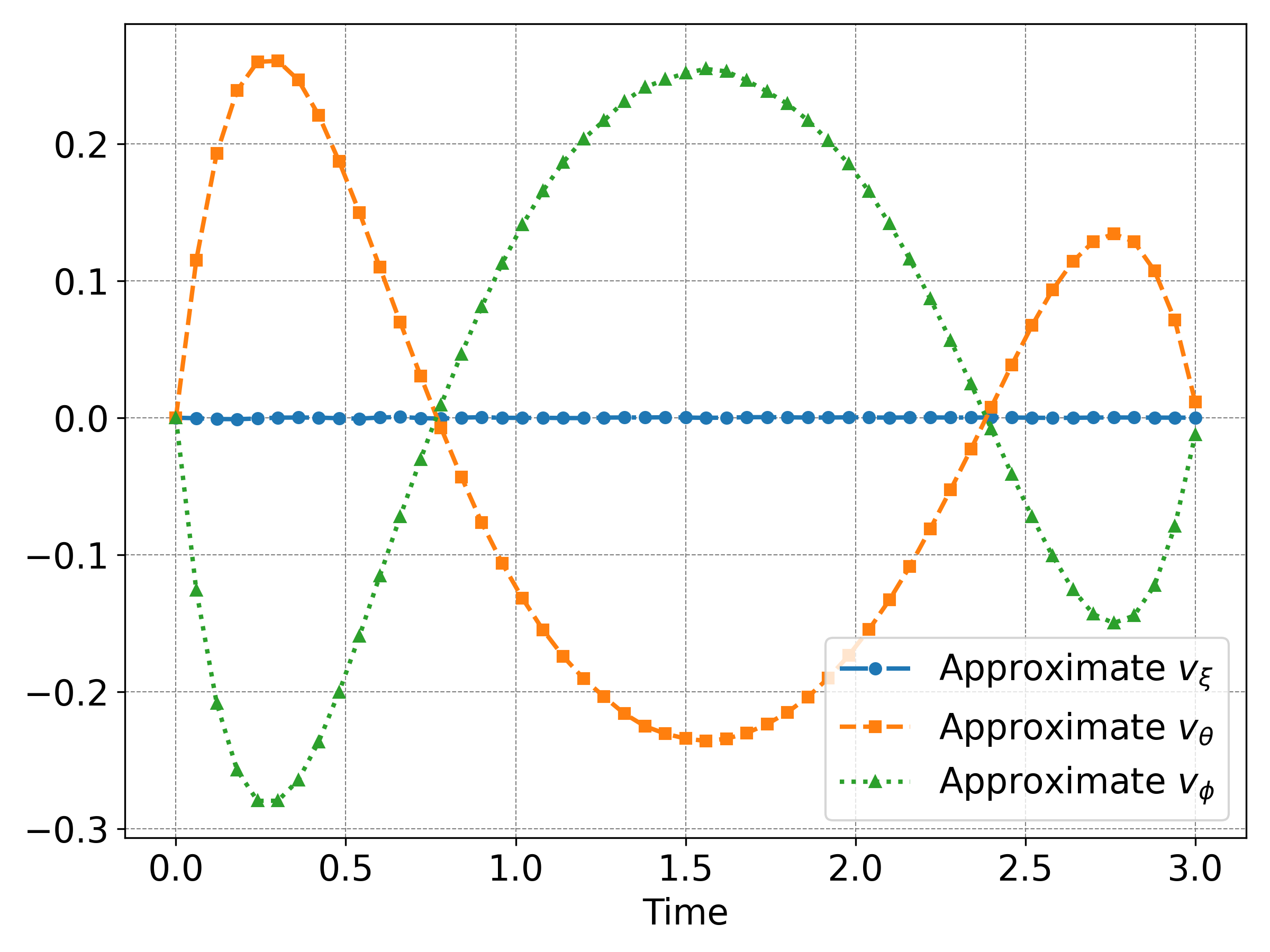}
			\caption{On the top we display the space coordinates of $50$ approximate optimal trajectories. The color blue represents the initial time, red indicates the final time, and intermediate colors correspond to intermediate times. Next, we display the profiles of the optimal controls and the remaining components of the state variable for one trajectory.} 
			\label{fig_quad}
		\end{figure}

		Our goal is to solve an optimal control problem to determine the optimal trajectories steering the quadcopter from an initial state towards a target state. We write system~\eqref{eq:quadcopter_dynamics} as a first-order controlled ODE by adding to the state of the system the velocity variables $v_{x}$, $v_{y}$, $v_{z}$, $v_{\xi}$, $v_{\theta}$, $v_{\phi}$, which yields a $12$-dimensional state space.  We consider as in~\cite{ lin2018splitting, lin2021alternating, onken2022neural}, running and terminal costs given respectively by 
		$$
		\ell(t,\boldsymbol{x},\alpha)=2+\left\|\alpha\right\|^2\quad\text{and}\quad g(\boldsymbol{x})=\frac{\|\boldsymbol{x}-\boldsymbol{\widehat{x}}\|^2}{\epsilon},
		$$
		where $\boldsymbol{x}=(x,y,z,v_{x},v_{y},v_{z},\xi,\theta,\phi,v_{\xi},v_{\theta},v_{\phi})$, $\alpha=(u,\tau_{\xi},\tau_{\theta},\tau_{\phi})$, $\epsilon>0$,  and $\boldsymbol{\widehat{x}}\in\RR^{12}$ is given. In our experiments, we set a time horizon $T=3$, a target point $\boldsymbol{\widehat{x}}=(1,1,1,0,0,0,0,0,0,0,0,0)$, and $\epsilon=0.01$. We implement the IGT method with a neural network $V_{\theta}^{{\rm NN}}$ parameterized with $\varphi_2$.  At the initialization step, the training of $V_{\theta}^{{\rm NN}}$ using the DGM with a batch sampled from $[0,T]\times\RR^{d}$ is computationally too expensive and hence we use the C-DGM variant to provide an initial estimate of the value function. Algorithm~\ref{alg:IGT} is run over two rounds with $\mathcal{D}_{\text{OC}}$ being generated with a batch of size $S=32$ of initial conditions and penalization parameter $\lambda_1=0.01$. The total training time is approximately 1464 seconds.  Figure~\ref{fig_quad} shows the approximate optimal trajectories for 
		problem $({\bf P}^{0,\boldsymbol{x}_{0}})$, where $\boldsymbol{x}_{0}=(x_{0},y_{0},z_{0},0,0,0,0,0,0,0,0,0)$ with $x_{0}$, $y_0$, $z_0$ being sampled randomly and independently from a Gaussian distribution of mean $-1$ and standard deviation $\sqrt{0.04}$. Let us point out that, contrary to the IGT method, in this example the time-marching technique considered in~\cite{kang2017mitigating,nakamura2021adaptive} to deal with the solution of the TPBVP problems failed to converge even when several  intermediate times are employed.
		
		\subsubsection{{\bf Obstacle problem}}  \label{subsec:obstacle_oc}
		We demonstrate the potential of the IGT strategy on a problem where the objective is to compute trajectories that drive agents from a given initial condition to a target state $\widehat{x} \in \mathbb{R}^d$, while avoiding obstacles and minimizing kinetic energy. Our goal is to show how successive rounds in the IGT algorithm can effectively refine convergence. To emphasize this, we formulate the problem to be highly sensitive to the initial guess, making the training process, particularly in the first step, potentially expensive, especially in high-dimensional cases.
		
		We consider the family of optimal control problems~\eqref{def:oc_problem} with $T=1$, $A=\mathbb{R}^{d}$, and  
		$$
		\ell(t,x,\alpha)= c \|\alpha\|^{2} + f(x),\quad 
		g(x)=\|x-\widehat{x}\|^2,\quad 
		b(t,x,\alpha)=-2 c \alpha\quad\text{for all }(t,x,\alpha)\in [0,T]\times\mathbb{R}^{d}\times\RR^{d},
		$$
		where $c>0$ and the function $f$ penalizes trajectories passing through restricted areas. The associated HJB equation with quadratic Hamiltonian reads:
		\begin{equation}
			\begin{aligned}
				-\partial_{t}V(t,x)+ c \left\|\nabla_{x}V(t,x)\right\|^2&= f(x)\quad\text{for all }(t,x)\in]0,1[\times\RR^{d}, \\
				V\left(1,x\right)&=g(x)\quad\text{for all }x\in\RR^{d}.
			\end{aligned}
			\label{hjb_oc_obs} 
		\end{equation}
		We consider circular obstacles in the first two coordinates of the state variable having the form 
		$$
		\mathcal{O}_{\mathbf{c},r}:=\{x\in\RR^{d}\,|\,f_{\mathbf{c},r}(x):=r^{2}-(x_{1}-\mathbf{c}_1)^{2}-(x_{2}-\mathbf{c}_2)^{2}\geq 0\},
		$$
		where $\mathbf{c}\in\RR^{2}$ and $r>0$. We penalize passing through a family of obstacles $\{\mathcal{O}_{\mathbf{c}^{i},r^{i}}\}_{i=1}^{N_{\text{o}}}$ by taking $f=\gamma_{\mathrm{obst}}\text{boltz}_{s}(f_{\mathbf{c}^1,r^1},\hdots,f_{\mathbf{c}^{N_{\text{o}}},r^{N_{\text{o}}}})$, where $\gamma_{\mathrm{obst}},\,s>0$. Here, $\text{boltz}_{s}(f_{\mathbf{c}^1,r^1},\hdots,f_{\mathbf{c}^{N_{\text{o}}},r^{N_{\text{o}}}})$ denotes the Boltzmann operator of the family $\{f_{\mathbf{c}^i,r^i}\}_{i=1}^{N_{\text{o}}}$, defined as 
		\begin{equation} \label{eq:boltz}
			\text{boltz}_{s}(f_{\mathbf{c}^1,r^1},\hdots,f_{\mathbf{c}^{N_{\text{o}}},r^{N_{\text{o}}}})(x)=\frac{\sum_{i=1}^{N_{\text{o}}}f_{\mathbf{c}^{i},r^{i}}(x)\exp\left(s f_{\mathbf{c}^{i},r^{i}}(x)\right)}{\sum_{i=1}^{N_{\text{o}}}\exp\left(s f_{\mathbf{c}^{i},r^{i}}(x)\right)}\quad\text{for all }x\in\RR^{d},
		\end{equation}
		which is a smooth function that approximates $\max_{i=,\hdots,N_{\text{o}}}f_{\mathbf{c}^i,r^i}$ for large $s$.
		
\begin{figure}
			\centering\includegraphics[width=5cm]{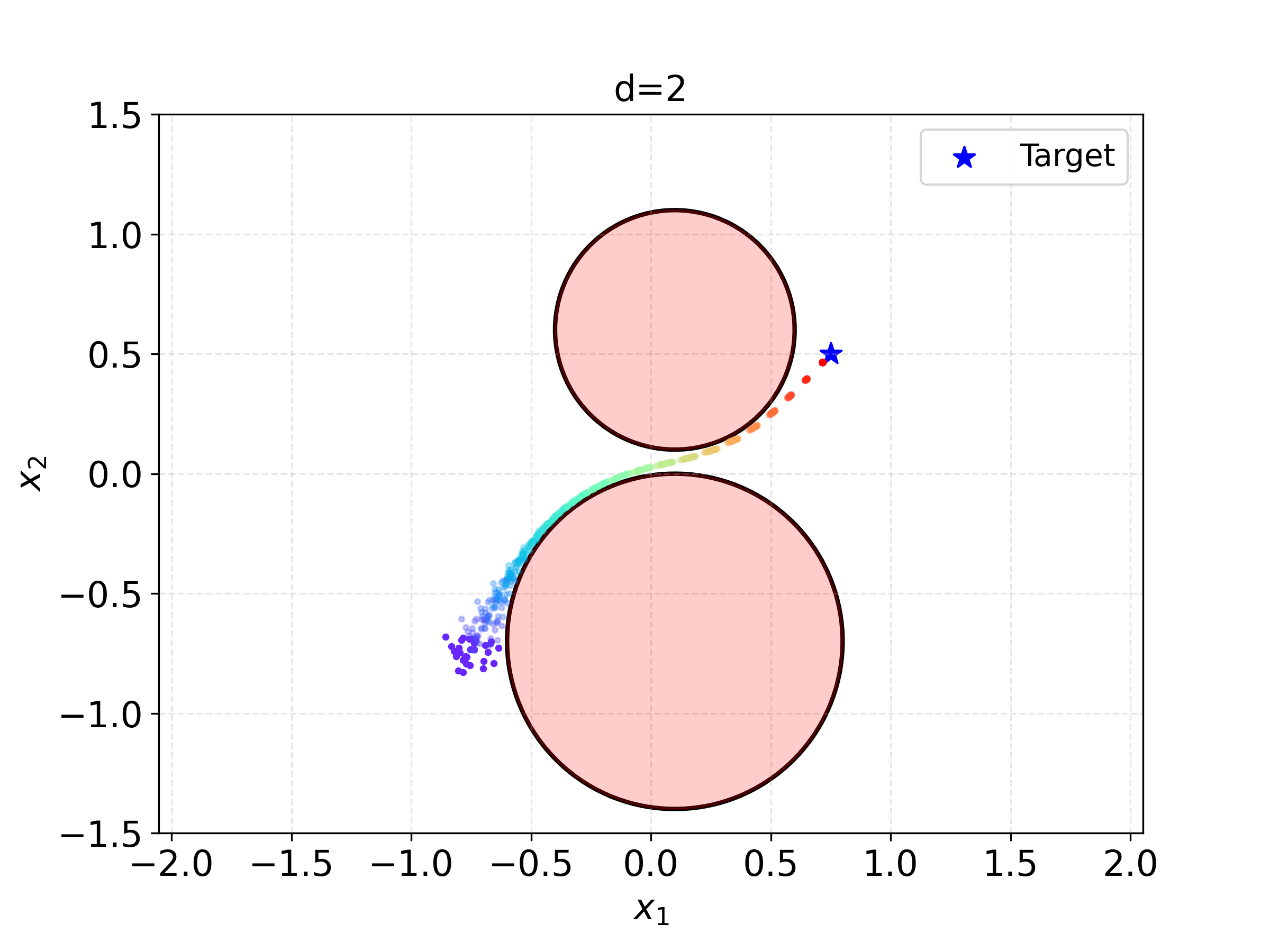}
			\centering\includegraphics[width=5cm]{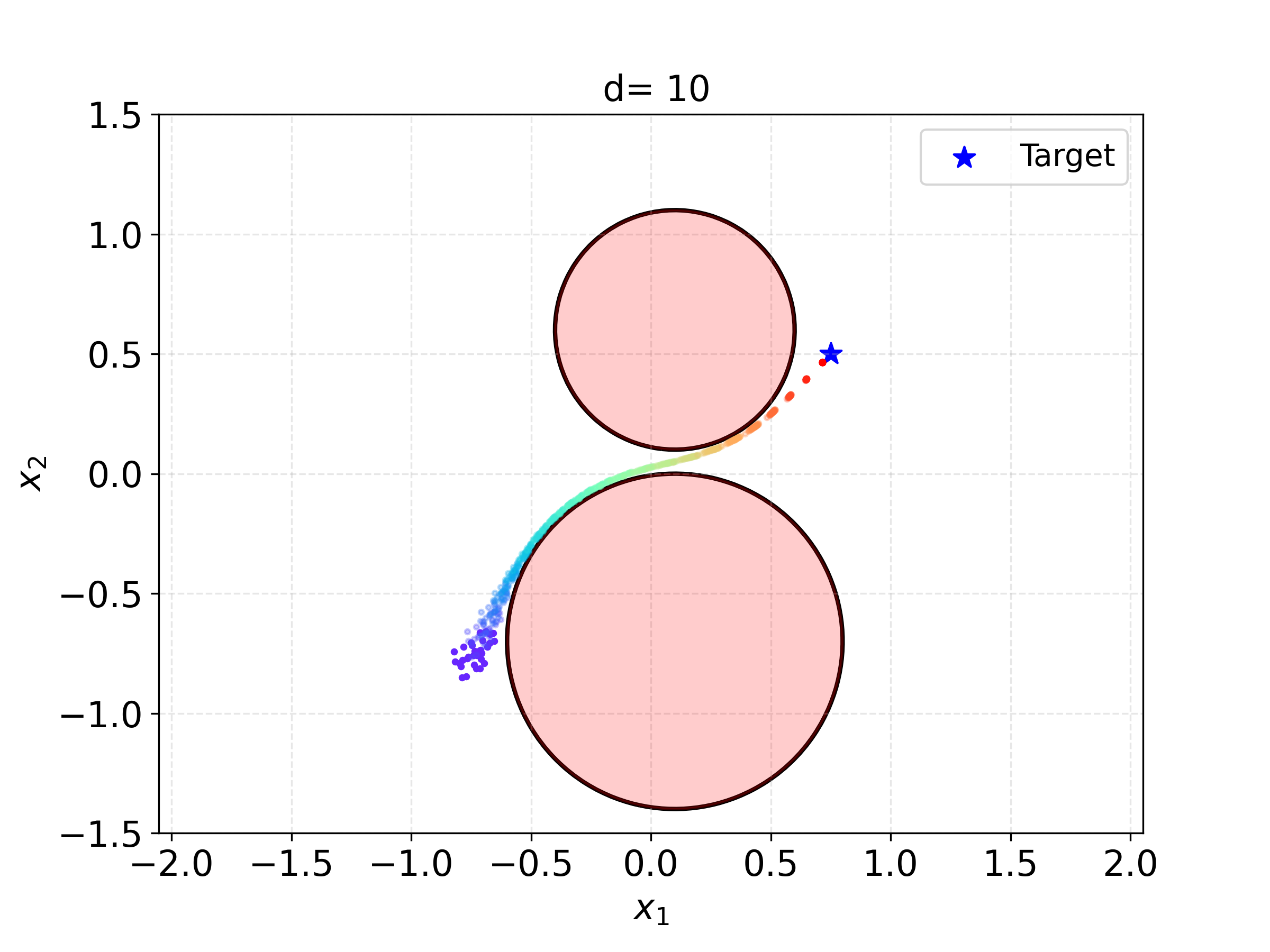}
			\centering\includegraphics[width=5cm]{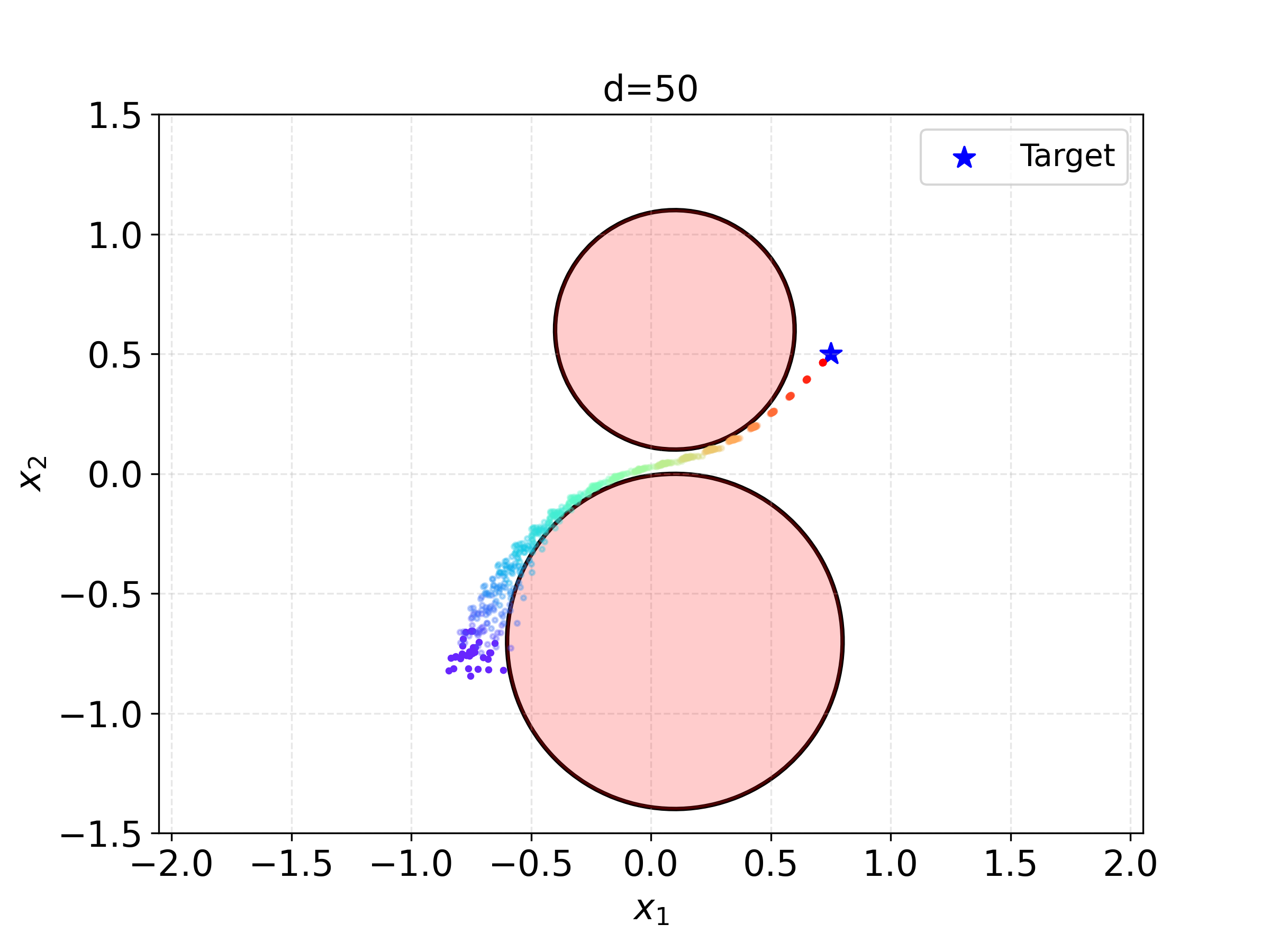}
			\caption{Projections on the first two coordinates of approximate optimal trajectories starting from $32$ initial conditions in dimensions $d=2,\,10,\,50$. Blue and red points indicate the first two coordinates at times $0$ and $1$, respectively.} 
			\label{fig: opt_obs}
		\end{figure}

		In our numerical tests, we take $N_{\text{o}}=2$, $\mathbf{c}^{1}=(0.1,0.6)$, $\mathbf{c}^{2}=(0,-0.7)$, $r^{1}=0.5$, $r^{2}=0.7$, $s=50$, $c=6$, $\gamma_{\mathrm{obst}}=5$, and target point $\boldsymbol{\widehat{x}}=(0.75,0.5,0,\ldots,0)$. We implement the IGT method with a neural network $V_{\theta}^{{\rm NN}}$ parameterized with $\varphi_2$.  We use the C-DGM at the initialization step, the dataset $\mathcal{D}_{\text{OC}}$ is generated with a batch of size $S=32$ of initial conditions and penalization parameter $\lambda_1=1$. As Table~\ref{tab:convergence_rounds} shows, to generate a reliable dataset $\mathcal{D}_{\text{OC}}$ we need to perform two rounds of the IGT method when $d=2,10$, and three rounds for $d=50$. This exemplifies the crucial role of the initialization step and the implementation of several rounds, as even if full convergence of the TPBVP fails in a subset of the batch of initial conditions, it is achieved by augmenting the number of rounds. On the other  hand, arbitrary initializations in solving the TPBVPs fail to converge at any of the points of the batch of initial conditions,  even when using time-marching with intermediate times, as it was the case in the previous example. Projections of the approximate optimal trajectories on the two first coordinates are displayed in Figure~\ref{fig: opt_obs}. We observe that the results are consistent across the three considered dimensions, confirming the robustness of the method with respect to the state dimension.
		
		\begin{table}[htbp]
			\caption{Convergence results for the TPBVP solver in the data generation step across dimensions $d=2,\,10,\,50$. Each entry shows the number of successful TPBVP  convergences out of 32 initial states.}
			\centering
			\begin{tabular}{|l|c|c|c|}
				\hline
				\text{Dimension} & \text{Round 1 (with initialization)} & \text{Round 2 (from Round 1)} & \text{Round 3 (from Round 2)} \\ \hline
				$d = 2$  & 14 / 32 & \text{32 / 32} & -- \\ \hline
				$d = 10$ & 8 / 32  & \text{32 / 32} & -- \\ \hline
				$d = 50$ & 3 / 32  & \text{21 / 32} & \text{32 / 32} \\ \hline
			\end{tabular}
			\label{tab:convergence_rounds}
		\end{table}
		
		\subsection{Mean field games}
		\label{subsec:numerical_mfg}
		
		In this section, we evaluate the effectiveness of the IGT-MFG algorithm by applying it to a series of examples. We begin with a problem for which the exact solution is known, allowing us not only to compute the exploitability and aggregated Sinkhorn divergences at each iteration, but also to measure the relative errors with respect to the exact solution.  
		We end this section by considering a MFG problem related to the obstacle optimal control problem discussed in the previous section. In some of these examples, we consider problems with different dimensions to highlight the robustness and effectiveness of the proposed method even in high-dimensional settings.
		
		In all cases, we consider the parametrization $\varphi_2$ for $V^{\rm NN}_{\theta}$ and $\Phi^{{\rm NN}}_{\omega}$. The maximum number of iterations of the fictitious play method $K_{\rm max}$ will vary depending on the type and complexity of the problem.
		\subsubsection{{\bf Evaluating IGT-MFG}} \label{subsec:ex_1_mfg}
		To evaluate the effectiveness of the IGT-MFG algorithm, we consider a linear quadratic MFG (see, e.g.,~\cite{MR3489817}) with an explicit solution which serves as a benchmark for numerical analysis. In this example, already considered in~\cite{Silva2024LG} to evaluate the performance of a Lagrange-Galerkin scheme for a first-order MFG, we look at system~\eqref{eq:MFG_system} with $\widetilde{\H}(x, p)=\|p\|^{2}/2$, coupling functions
		$$
		F(x,m)=\frac{1}{2}\left\|x-\int_{\mathbb{R}^d} y \mathrm{~d} m(y)\right\|^2, \quad G(x,m)=0,
		$$
		and initial data $m_0$ given by a $d$-dimensional Gaussian distribution with mean $\mathbf{a}\in\mathbb{R}^d$ and covariance matrix $\Sigma_0\in\mathbb{R}^{d \times d}$. For simplicity we assume that $\Sigma_0$ is diagonal. Setting $I_{d}$ for the $d\times d$ identity matrix and
		$$
		\Pi(t)=\left(\frac{e^{2 T-t}-e^t}{e^{2 T-t}+e^t}\right) I_d, \quad s(t)=-\Pi(t) \mathbf{a}, \quad c(t)=\frac{1}{2} (\Pi(t) \mathbf{a}) \cdot \mathbf{a} \quad \text{ for all } t \in[0, T],
		$$
		system~\eqref{eq:MFG_system} admits a unique solution $\left(V^*, m^*\right)$, where 
		$$
		V^*(t, x)=\frac{1}{2} (\Pi(t)x) \cdot x +  s(t) \cdot x + c(t) \quad \text { for all }(t, x) \in[0, T] \times \mathbb{R}^d
		$$
		and, for every $t \in[0, T]$, $m^*(t)$ is a $d$-dimensional Gaussian distribution with mean $\mathbf{a}$ and diagonal covariance matrix given by 
		$$
		\Sigma_t=\left(\frac{e^{2 T-t}+e^t}{e^{2 T}+1}\right)^2\Sigma_0.
		$$
		
		Since the exact solution of this problem is known, our goal in these experiments is to illustrate how the relative errors, exploitability, and agregated Sinkhorn divergences evolve over the iterations of the algorithm. In the numerical tests below we set, for every $l=1,\hdots,d$, $\mathbf{a}_l=0.1$, $\left(\Sigma_0\right)_{l,l}=0.105$, and $T=1$. We test the method in the one dimensional case $d=1$  and also when $d=10, 50,$ and  $100$. For $d=1$, we take a batch of initial conditions $\B$ of size $B=1000$ sampled from $m_0$, while for $d=10$, we use a larger batch of size $B=4000$ and $B=10000$ for $d=50, 100$.  This batch size has been selected larger in order to better approximate the coupling term $F$.

		{\it One  dimensional test.} In this test, we consider the case $d=1$ and initialize the IGT-MFG method by performing a joint training of the value function and the flow map. This initialization employs a batch of size $M_1'=256$, sampled uniformly from the time-space domain $[0,1] \times [-2,2]$, and a batch of trajectories of size $M_2'=256$ for the penalization term $\lossHJB^0$ in~\eqref{eq:loss_hjb_mfg_init}. Additionally, for the penalization term $\lossODE^0$ in~\eqref{eq:loss_ode_mfg_init}, we use a batch of size $M_3'=256$ sampled from the initial distribution $m_0$. At each iteration $k$ of the fictitious play method, we update $\theta^*$ using the IGT algorithm (see Algorithm~\ref{alg:IGT}). In the initialization step of the IGT algorithm, we consider a batch of uniformly distributed points in the time-space domain $[0,1] \times [-2,2]$ of size $M_1=256$. In the data generation step,  we take a batch of size $S_1=128$  of initial conditions sampled from $m_0$. In the training step, the computation of the penalization term $\lossHJB$ in~\eqref{eq:loss_total} is performed using the same size $M_1$ of uniformly distributed points in $[0,1] \times [-2,2]$ and penalization parameter $\lambda_1 = 1$.  Then to compute $\Phi^{\rm NN}_{\omega^*}$ we first use a batch of size $S_2=256$, sampled from $m_0$, to generate data. In the training step, we use a batch of size $M_2=256$ to compute $\lossODE$ and we take $\lambda_2 = 0.5$ in the definition of~$\mathrm{Loss}_{\mathrm{gen}}$ (see~\eqref{eq:loss_total_gen}).
		
		Table~\ref{tab: MFG1}  provides the relative errors $E_{\infty}(V^{\mathrm{NN}}_{\theta}(t,\cdot))$ and $E_{2}(V^{\mathrm{NN}}_{\theta}(t,\cdot))$, along with exploitability $\psi^{\B}(\alpha_{\theta^*})$ (see~\eqref{eq:exploitability_numerical}) and aggregate Sinkhorn divergences (see~\eqref{eq:w_metrics}), with $\varepsilon=1/2$, between $\overline{m}_k$ and its best response $\mu_{k+1}$, between $\overline{m}_k$ and $m^*$, and between the best response $\mu_{k+1}$ and $m^*$. 
		We notice that we have needed $k=8$ iterations to achieve the desired tolerance $\texttt{tol}=5\times10^{-6}$ for ${\bf S_{\varepsilon}^{\infty}}(\overline{m}^k,\mu_{k+1})$.

		\begin{table}[htbp]
			\caption{Relative errors of the value function, exploitability, and aggregated Sinkhorn divergences at each fictitious play iteration $k$ for $d=1$.}
			\renewcommand{\arraystretch}{1.05}
			\setlength{\tabcolsep}{4pt}
			\centering
			\begin{tabular}{|c|c|c|c|c|c|c|}
				\hline
				\multirow{2}{*}{$k$} 
				& \multicolumn{2}{c|}{Relative error $E_\infty$}
				&  \multirow{2}{*}{$\psi^{\mathcal{B}}(\alpha_{\theta^*})$}
				&  \multirow{2}{*}{${\bf S_{\varepsilon}^\infty}(\overline{m}^k,\mu_{k+1})$} 
				&  \multirow{2}{*}{${\bf S_{\varepsilon}^\infty}(\overline{m}^k,m^*) $}
				&  \multirow{2}{*}{${\bf S_{\varepsilon}^\infty}(\mu_{k+1},m^*) $} \\
				\cmidrule{2-3}
				& $t=0$ & $t=0.5$ & & & & \\
				\hline
				1 & $2.996{\times}10^{-3}$ & $4.515{\times}10^{-3}$ & $4.50{\times}10^{-2}$ & $6.65{\times}10^{-5}$ & $2.30{\times}10^{-3}$ & $1.50{\times}10^{-3}$ \\
				2 & $8.214{\times}10^{-4}$ & $2.392{\times}10^{-3}$ & $4.98{\times}10^{-3}$ & $2.60{\times}10^{-5}$ & $1.00{\times}10^{-3}$ & $7.86{\times}10^{-4}$ \\
				3 & $4.904{\times}10^{-4}$ & $1.576{\times}10^{-3}$ & $4.32{\times}10^{-3}$ & $1.30{\times}10^{-5}$ & $1.12{\times}10^{-3}$ & $7.69{\times}10^{-4}$ \\
				4 & $7.877{\times}10^{-4}$ & $1.144{\times}10^{-3}$ & $3.61{\times}10^{-3}$ & $7.38{\times}10^{-6}$ & $1.04{\times}10^{-3}$ & $5.96{\times}10^{-4}$ \\
				5 & $3.516{\times}10^{-4}$ & $1.411{\times}10^{-3}$ & $3.89{\times}10^{-3}$ & $5.56{\times}10^{-6}$ & $9.42{\times}10^{-4}$ & $9.29{\times}10^{-4}$ \\
				6 & $5.833{\times}10^{-4}$ & $7.160{\times}10^{-4}$ & $2.55{\times}10^{-3}$ & $5.78{\times}10^{-6}$ & $8.77{\times}10^{-4}$ & $7.82{\times}10^{-4}$ \\
				7 & $6.786{\times}10^{-4}$ & $8.433{\times}10^{-4}$ & $1.57{\times}10^{-3}$ & $5.28{\times}10^{-6}$ & $8.65{\times}10^{-4}$ & $6.37{\times}10^{-4}$ \\
				8 & $4.906{\times}10^{-4}$ & $7.358{\times}10^{-4}$ & $1.05{\times}10^{-3}$ & $3.32{\times}10^{-6}$ & $6.27{\times}10^{-4}$ & $7.85{\times}10^{-4}$ \\
				\hline
			\end{tabular}
			\label{tab: MFG1}
		\end{table}

		{\it High-dimensional test.} We now examine the same experiment in higher dimensions $d=10, 50, 100$. We initialize the IGT-MFG method using batch sizes $M_1'=M_2'=M_3'=512$. Throughout the algorithm, we maintain batch sizes $M_1$, $M_2$, and $S_2$ to $512$. For the data generation step in the IGT method, we set $S_1=128$ for $d=10$, and increase it to $S_1=512$ for $d=50$ and $d=100$. Regarding the penalization parameter, we choose $\lambda_1 = 0.5$ for $d=10$, $\lambda_1 = 0.1$ for $d=50$, and $\lambda_1 = 0.01$ for $d=100$. All uniform samplings are taken in the time-space domain $[0,1] \times [-2,2]^d$.  Table~\ref{tab:multi_dim_comparison} summarize the performance and convergence behavior of the algorithm with $K_{max}=20$ fictitious play iterations.

		\begin{table}[htbp]
			\caption{Relative errors of the value function, exploitability, and aggregated Sinkhorn divergences at selected fictitious play iterations $k$ for dimensions $d=10,50,100$.}
			\renewcommand{\arraystretch}{1.1}
			\setlength{\tabcolsep}{4pt}
			\centering
			\begin{tabular}{|c|c|c|c|c|c|c|}
				\hline
				$d$ & $k$ 
				& $E_\infty(t=0)$ 
				& $\psi^{\mathcal{B}}(\alpha_{\theta^*})$
				& ${\bf S_{\varepsilon}^\infty}(\overline{m}^k,\mu_{k+1})$
				& ${\bf S_{\varepsilon}^\infty}(\overline{m}^k,m^*)$
				& ${\bf S_{\varepsilon}^\infty}(\mu_{k+1},m^*)$ \\
				\hline
				\multirow{6}{*}{10}
				& 1  & $4.66{\times}10^{-2}$ & $8.28{\times}10^{-2}$ & $2.37{\times}10^{-3}$ & $2.50{\times}10^{-1}$ & $2.43{\times}10^{-1}$ \\
				& 4  & $5.07{\times}10^{-2}$ & $9.72{\times}10^{-2}$ & $5.33{\times}10^{-4}$ & $2.48{\times}10^{-1}$ & $2.42{\times}10^{-1}$ \\
				& 8  & $2.33{\times}10^{-2}$ & $5.01{\times}10^{-2}$ & $2.94{\times}10^{-4}$ & $2.45{\times}10^{-1}$ & $2.43{\times}10^{-1}$ \\
				& 12 & $2.87{\times}10^{-2}$ & $2.07{\times}10^{-2}$ & $2.23{\times}10^{-4}$ & $2.40{\times}10^{-1}$ & $2.47{\times}10^{-1}$ \\
				& 16 & $2.44{\times}10^{-2}$ & $3.61{\times}10^{-2}$ & $1.73{\times}10^{-4}$ & $2.41{\times}10^{-1}$ & $2.48{\times}10^{-1}$ \\
				& 20 & $2.05{\times}10^{-2}$ & $3.48{\times}10^{-2}$ & $1.32{\times}10^{-4}$ & $2.42{\times}10^{-1}$ & $2.40{\times}10^{-1}$ \\
				\hline
				\multirow{6}{*}{50}
				& 1  & $9.00{\times}10^{-1}$ & $2.21{\times}10^{1}$ & $2.46{\times}10^{-1}$ & $8.20{\times}10^{-1}$ & $7.87{\times}10^{-1}$ \\
				& 4  & $3.24{\times}10^{-1}$ & $2.54{\times}10^{0}$ & $6.92{\times}10^{-2}$ & $7.68{\times}10^{-1}$ & $7.76{\times}10^{-1}$ \\
				& 8  & $1.40{\times}10^{-1}$ & $6.12{\times}10^{-1}$ & $3.15{\times}10^{-2}$ & $7.79{\times}10^{-1}$ & $7.81{\times}10^{-1}$ \\
				& 12 & $8.96{\times}10^{-2}$ & $1.11{\times}10^{0}$ & $2.53{\times}10^{-2}$ & $7.77{\times}10^{-1}$ & $7.68{\times}10^{-1}$ \\
				& 16 & $5.95{\times}10^{-2}$ & $3.22{\times}10^{-1}$ & $1.85{\times}10^{-2}$ & $7.82{\times}10^{-1}$ & $7.75{\times}10^{-1}$ \\
				& 20 & $1.25{\times}10^{-1}$ & $2.73{\times}10^{-1}$ & $1.58{\times}10^{-2}$ & $7.71{\times}10^{-1}$ & $7.70{\times}10^{-1}$ \\
				\hline
				\multirow{6}{*}{100}
				& 1  & $4.87{\times}10^{-1}$ & $2.27{\times}10^{1}$ & $1.41{\times}10^{0}$ & $1.59{\times}10^{0}$ & $2.35{\times}10^{0}$ \\
				& 4  & $2.07{\times}10^{-1}$ & $1.48{\times}10^{1}$ & $4.19{\times}10^{-1}$ & $1.55{\times}10^{0}$ & $1.53{\times}10^{0}$ \\
				& 8  & $1.12{\times}10^{-1}$ & $5.56{\times}10^{0}$ & $1.98{\times}10^{-1}$ & $1.53{\times}10^{0}$ & $1.54{\times}10^{0}$ \\
				& 12 & $1.16{\times}10^{-1}$ & $3.35{\times}10^{0}$ & $1.30{\times}10^{-1}$ & $1.53{\times}10^{0}$ & $1.55{\times}10^{0}$ \\
				& 16 & $9.92{\times}10^{-2}$ & $2.56{\times}10^{0}$ & $1.14{\times}10^{-1}$ & $1.54{\times}10^{0}$ & $1.54{\times}10^{0}$ \\
				& 20 & $8.78{\times}10^{-2}$ & $2.45{\times}10^{0}$ & $1.01{\times}10^{-1}$ & $1.53{\times}10^{0}$ & $1.54{\times}10^{0}$ \\
				\hline
			\end{tabular}
			\label{tab:multi_dim_comparison}
		\end{table}

		\subsubsection{{\bf Example with obstacles}}
		We illustrate the potential of the IGT-MFG method on a MFG in which, for a population having aversion to crowded places, the objective is to obtain trajectories that steer agents from a given initial distribution to a target state $ \widehat{x} \in \mathbb{R}^d $, while avoiding obstacles and minimizing kinetic energy. Although this setting shares structural similarities with the optimal control example previously discussed in Subsection \ref{subsec:obstacle_oc}, particularly in terms of the dynamics and target-seeking behavior, it fundamentally differs due to the inclusion of the aversion to crowd effect. The latter introduces interactions among agents via their cost functions, leading to a more complex and realistic modeling framework. In this context, based on the results presented in Subsection \ref{subsec:obstacle_oc}, we observe that successive rounds of the IGT algorithm improve convergence and yield a more accurate approximation of the value function.
		
		We study the MFG system~\eqref{eq:MFG_system} with quadratic Hamiltonian $\widetilde{\mathcal{H}}(x,p)=c\,\|p\|^{2}$
		and coupling functions
		$$
		F(x,m) \;=\; \theta_{1}\,f(x) \;+\; \theta_{2}\,\bigl(\rho\ast\tilde{m}\bigr)(x_{1},x_{2}),
		\qquad
		G(x,m) \;=\; \frac{\theta_{3}}{2}\,\|x-\widehat{x}\|^{2},
		$$
		where the weights satisfy \(\theta_{1},\theta_{2},\theta_{3}\ge 0\), $d\geq 2$ and $x_1,\, x_2$ are the first two coordinates of $x\in\RR^d$.  
		The function $f\colon\mathbb{R}^{d}\!\to\mathbb{R}$ penalises trajectories that enter restricted regions, while the aversion to crowd effect is captured by the convolution term. Here $m\in\P(\RR^d)$, $\tilde{m}$ denotes its marginal distribution  with respect to its first $2$ variables,  and $\rho$ is the density of a $2$-dimensional Gaussian vector with $0$ mean and covariance matrix $\Sigma\in\RR^{2\times 2}$ . We restrict the motion in the first two spatial coordinates by prescribing a collection of smooth functions
		$$
		f_{\mathbf{c}_{i},R_{i},Q_{i},b_{i}}(x)\;:=\;
		-\;v_{i}^{\top}Q_{i}v_{i}\;-\;b_{i}\!\cdot\!v_{i}\;-\;1,
		\qquad 
		v_{i}\;=\;\bigl((x_{1},x_{2})-\mathbf{c}_{i}\bigr)R_{i},
		\quad i=1,\dots,N_{\mathrm{o}},
		$$
		where, $\mathbf{c}_{i}\in\RR^{2}$ is the centre of the $i$-th obstacle,  $R_{i}\in\mathrm{SO}(2)$ is a rotation matrix (e.g.\ $R_{i}=R(\theta_{i})$ with angle $\theta_{i})$, $Q_{i}\in\mathbb{S}^{2}_{+}$ is a positive semi-definite shape matrix, and  
		$b_{i}\in\RR^{2}$ tilts the quadratic form to create asymmetric profiles. The corresponding obstacle region is the super-level set  
		$$
		\mathcal{O}_{i}\;=\;\bigl\{x\in\RR^{d}\;:\;f_{\mathbf{c}_{i},R_{i},Q_{i},b_{i}}(x)\;\ge 0\bigr\},
		\qquad i=1,\dots,N_{\mathrm{o}}.
		$$ We penalize passing through a family of obstacles $\{\mathcal{O}_i\}_{i=1}^{N_{\text{o}}}$ by taking 
		$$
		f(x)\;=\;
		\mathrm{boltz}_{s}\!\bigl(f_{\mathbf{c}_{1},R_{1},Q_{1},b_{1}}(x),\dots,
		f_{\mathbf{c}_{N_{\mathrm{o}},R_{N_{\mathrm{o}}},Q_{N_{\mathrm{o}}},b_{N_{\mathrm{o}}}}}(x)\bigr),
		\qquad s>0,
		$$ where, $\text{boltz}_{s}$ denotes the Boltzmann operator as in \eqref{eq:boltz}.
		
		In our numerical tests, we take $T=1$, $c=3/2$, $\theta_{1}=\theta_{2}=\theta_{3}=3$, $\Sigma = (0.25) I_{2}$, and target point ${\widehat{x}}=(0.75,0,0,\ldots,0)$. The parameters involved in the Boltzmann operator are  $N_{\text{o}}=2$, $\mathbf{c}_{1}=(0,0.3)$, $\mathbf{c}_{2}=(0,-0.3)$, $R_{1}=R_2=\left(\begin{array}{rr}\cos (\theta) & -\sin (\theta) \\ \sin (\theta) & \cos (\theta)\end{array}\right)$ with $\theta=\pi / 5,$  $Q_1=Q_2=\left(\begin{array}{ll}10 & 0 \\ 0 & 1\end{array}\right)$, $b_1=(0,3)$, $b_2=(0, -3)$, and  $s=50$. Our initial density $m_0$ is a Gaussian centered at $(-0.75,0,0,\ldots,0)\in \RR^d$ with covariance matrix $\Sigma_0=(0.1) I_d$. 
		
		We test the method in the two-dimensional case $d = 2$ and also when $d = 10$.  As mentioned above, taking into account the results of Subsection \ref{subsec:obstacle_oc}, multiple rounds of IGT are typically required. Due to the computational cost of this, we have considered batches $\B$ of smaller sizes than those used, for example, in Example \ref{subsec:ex_1_mfg}. For $d = 2$, we use a  batch $\B$ of size $B = 264$ sampled from $m_0$, while for $d = 10$, we employ a batch of size $B = 512$. In both cases, we have used a tolerance $\texttt{tol}=10^{-2}$ as the stopping criterion for Algorithm \ref{alg: IGT-MFG}.

\begin{figure}[h]
			\centering\includegraphics[width=6cm]{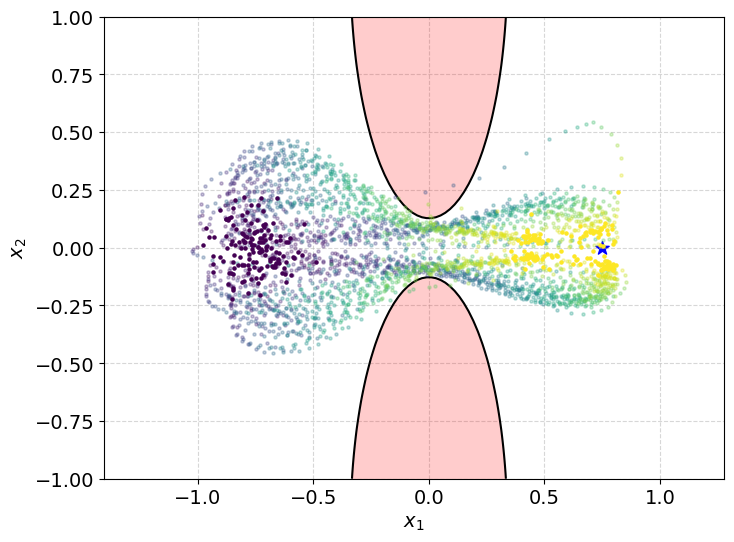}
			\centering\includegraphics[width=6cm]{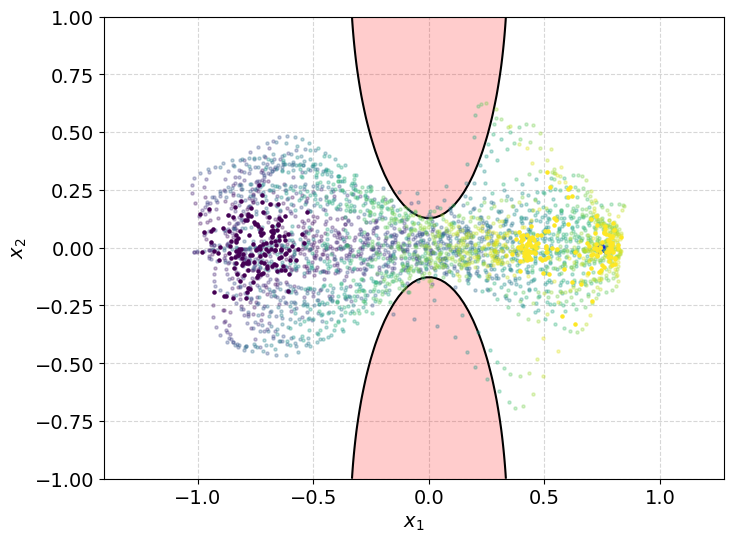}
			\caption{Projection onto the first two coordinates of approximate optimal trajectories at equilibrium starting from $150$ initial conditions for $d = 2$ (left) and $d = 10$ (right). Blue and yellow points indicate positions at times $t = 0$ and $t = 1$, respectively
			} 
			\label{fig: MFG_obs}
		\end{figure}

		{\it Two  dimensional test.}  In this test, we consider the case $d=2$ and initialize the IGT-MFG method by performing a joint training of the value function and the flow map. This initialization employs a batch of size $M_1'=128$, sampled from the initial distribution $m_0$ for the penalization term $\lossODE^0$ in~\eqref{eq:loss_ode_mfg_init}, and a batch of size $M_2'=128$ sampled uniformly from the time-space domain $[0,1] \times [-1,1]^2$ for the penalization term $\lossHJB^0$ in~\eqref{eq:loss_hjb_mfg_init}. At each iteration $k$ of the fictitious play method, we update $\theta^*$ using the IGT algorithm (see Algorithm~\ref{alg:IGT}). The DGM initialization employs a batch of $M_{1} = 128$ uniformly distributed points in the time-space domain $[0,1] \times [-1,1]^2$. We then sample a batch of $S_{1} = 32$ initial conditions from $m_{0}$ to form the dataset $\mathcal{D}_{\mathrm{OC}}$. The penalization term $\lossHJB$ in~\eqref{eq:loss_total} is evaluated with the same size $M_1$ of points uniformly distributed in the same domain and weight $\lambda_{1} = 1$. To compute the generator $\Phi^{\mathrm{NN}}_{\omega^{\ast}}$ we use a batch of size $S_{2} = 128$ from $m_{0}$ to build $\mathcal{D}_{\mathrm{MFG}}$, together with $M_{2} = 128$ points for the penalization term $\lossODE$ in~\eqref{eq:loss_total_gen}, employing the weight $\lambda_{2} = 0.01$. On the left of Figure~\ref{fig: MFG_obs} we see the approximated optimal trajectories. Table~\ref{tab:obs_2} provides the exploitability and the aggregated Sinkhorn divergences across the iterations.

		\begin{table}[htbp]
			\caption{Exploitability and aggregated Sinkhorn divergence ${\bf S^{\infty}_{\varepsilon}}(\overline{m}_k,\mu_{k+1})$ at selected fictitious play iterations $k$ for $d=2$.}
			\setlength{\tabcolsep}{4pt}
			\renewcommand{\arraystretch}{1.2}
			\centering
			\begin{tabular}{|c|c|c|}
				\hline
				$k$ 
				& $\psi^{\mathcal{B}}(\alpha_{\theta^*})$ 
				& ${\bf S^{\infty}_{\varepsilon}}(\overline{m}^k,\mu_{k+1})$ \\
				\hline
				1  & $3.43574\times10^{-1}$ & $4.136\times10^{-1}$ \\
				5  & $1.30008\times10^{-1}$ & $1.154\times10^{-1}$ \\
				10  & $3.0416\times10^{-2}$  & $4.516\times10^{-2}$ \\
				15 & $3.6668\times10^{-2}$  & $2.597\times10^{-2}$ \\
				20 & $2.0108\times10^{-2}$  & $1.754\times10^{-2}$ \\
				\hline
			\end{tabular}
			\label{tab:obs_2}
		\end{table}

		{\it High-dimensional test.} We now consider the same experiment in a higher dimension $d=10$. Given the increased complexity of high-dimensional settings, we double the batch sizes to $M_1'=M_2'=M_3'=256$, and $M_1 = 256$, $S_1 = 64$, $S_2 = 256$, and $M_2 = 256$. Furthermore, we set the penalization parameters to $\lambda_1=0.5$ and $\lambda_2=0.001$. On the right of Figure~\ref{fig: MFG_obs} we display the projection   onto the first two coordinates of approximate optimal trajectories at equilibrium. Table~\ref{tab: obs_10} shows the exploitability and the  aggregated Sinkhorn divergences for this example.
		
		\begin{table}[htbp]
			\caption{Exploitability and aggregated Sinkhorn divergence ${\bf S^{\infty}_{\varepsilon}}(\overline{m}_k,\mu_{k+1})$ at selected fictitious play iterations $k$ for $d=10$.}
			\setlength{\tabcolsep}{4pt}
			\renewcommand{\arraystretch}{1.2}
			\centering
			\begin{tabular}{|c|c|c|}
				\hline
				$k$ 
				& $\psi^{\mathcal{B}}(\alpha_{\theta^*})$ 
				& ${\bf S^{\infty}_{\varepsilon}}(\overline{m}^k,\mu_{k+1})$ \\
				\hline
				1  & $5.90511\times10^{-1}$ & $4.892\times10^{-1}$ \\
				5  & $3.17829\times10^{-1}$ & $1.394\times10^{-1}$ \\
				10  & $1.85393\times10^{-1}$ & $3.896\times10^{-2}$ \\
				15 & $2.14395\times10^{-1}$ & $2.387\times10^{-2}$ \\
				20 & $8.3636\times10^{-2}$  & $1.801\times10^{-2}$ \\
				\hline
			\end{tabular}
			\label{tab: obs_10}
		\end{table}
		
		\FloatBarrier
		\section*{Acknowledgments}
		JG was partially supported by project MATH AmSud 23-MATH-17 (VIPS) and by PICT-2021-I-INVI-00834 (ANPCyT). FJS was partially supported by l'Agence Nationale de la Recherche (ANR), project ANR-22-CE40-0010, by KAUST through the subaward agreement ORA-2021-CRG10-4674.6, and by Minist\`ere de l'Europe et des Affaires \'etrang\`eres (MEAE), project MATH AmSud 23-MATH-17.
		
		\bibliographystyle{IGT}
		\bibliography{references}
		
	\end{document}